\documentclass[12pt]{article}
\usepackage{graphicx}
\usepackage{amsfonts}
\usepackage{subfigure}
\usepackage{amsmath,amssymb}

\renewcommand{\vec}[1]{\mbox{\boldmath$#1$}}

\DeclareMathOperator{\sdiv}{div}

\DeclareMathOperator{\cof}{cof}

\newtheorem{theorem}{Theorem}[section]

\newtheorem{lemma}[theorem]{Lemma}
\newtheorem{remark}[theorem]{Remark}

\usepackage{hyperref}
\hypersetup{bookmarks=true, bookmarksnumbered=true,
bookmarksopen=true, colorlinks, linkcolor=blue, citecolor=blue,
plainpages=false, pdfstartview=FitH}

\numberwithin{equation}{section}

\begin{document}
\pagestyle{plain}

\title{ A Meshing Strategy for a Quadratic Iso-parametric
FEM in Cavitation Computation in Nonlinear Elasticity\thanks{The research was
supported by the NSFC
project 11171008 and RFDP of China.}}

\author{Chunmei Su \hspace{1mm} and \hspace{1mm} Zhiping Li\thanks{Corresponding author,
email: lizp@math.pku.edu.cn} \\ \hspace{2mm} \\ LMAM \& School of Mathematical Sciences, \\
Peking University, Beijing 100871, China}

\date{}
\maketitle
\begin{abstract}
The approximation properties of a quadratic iso-parametric
finite element method for a typical cavitation problem in nonlinear elasticity
are analyzed. More precisely, (1) the finite element interpolation errors are
established in terms of the mesh parameters; (2) a mesh distribution strategy based on
an error equi-distribution principle is given; (3) the convergence of finite element
cavity solutions is proved. Numerical experiments show that, in fact, the optimal
convergence rate can be achieved by the numerical cavity solutions.
\end{abstract}

{\small \noindent

}

%\newpage
%\tableofcontents
%\newpage

\section{Introduction}
Nonlinear soft elastic materials, such as polymers, biological tissues, rubbers, etc.,
can display a particular singular deformation, referred to in the literature as cavitation,
when strong external force is applied \cite{Cristiano, Dorfmann, Gent and Lindley,
Kundu, Michel}. The occurrence and growth of cavities is considered closely related to the
material instability and to the damage and failure mechanisms of the materials
\cite{Adel,  Xu2015, Xu, Jaravel, Victor}. A huge number of work has been done by
numerous authors analysing cavitation experimentally, theoretically as well as numerically.

Generally speaking, there are two representative approaches characterizing cavitation.
One is the so-called defect model, which is based on the hypothesis that cavities grow
from pre-existing micro defects. Under this assumption, Gent and Lindley
\cite{Gent and Lindley} analysed the critical hydrostatic pressure at which
a given unit spherical void in an infinite extension of a Neo-Hookean material
would blow up, which was in a good agreement with their experiments therein. The other
is the perfect model established by Ball \cite{Ball82} based on the analytical evidence
that, under certain circumstances, a deformation with cavities created in an
originally intact material can be energetically favorable. It is shown that, under the
assumption that the cavities can appear only at a finite number of fixed points in the
intact materials, the solution of the defect model converges to the solution of the
perfect model \cite{Henao2009, Sival2006} as the radii of the pre-existing small voids
go to zero. In addition, analytical and numerical evidences indicate that whether a
point can serve as a possible position of cavitation can be evaluated by calculating
the corresponding configurational forces \cite{Lian and Li iso, Sival2002}.

The perfect model typically exhibits the Lavrentiev phenomenon \cite{Lavrentiev} when
there is a cavitation solution, leading to the failure of the conventional finite element
methods \cite{Bai and Li, Ball and Knowles}.
Though there are existing numerical methods developed to deal with the Lavrentiev
phenomenon \cite{Bai and Li, Ball and Knowles, Z. Li, Marrero and Betancourt},
they do not seem to be suitable to tackle the cavitation problem.
In fact, most of the numerical studies on cavitation are based on the defect model,
in which one considers to minimize the total energy of the form
\begin{equation}
  \label{functional}
  E(\vec{u})=\int_{\Omega_\varrho}W(\nabla \vec{u}(\vec{x}))d\vec{x},
\end{equation}
in the set of admissible functions
\begin{equation}
  \label{admissible set}
  \mathcal{A}=\{\vec{u}\in W^{1,p}(\Omega_\varrho;\mathbb{R}^n) \ \mbox{is one-to-one
  a.e.}: \vec{u}|_{\Gamma_0}=\vec{u}_0, \det \nabla \vec{u} >0 \ a.e. \},
\end{equation}
where $\Omega_\varrho = \Omega \setminus \bigcup_{i=1}^{m}
B_{\varrho_i}(\vec{a}_i)\subset\mathbb{R}^{n}\,(n=2,3)$
denotes the region occupied by an elastic body in its reference configuration,
$B_{\varrho_i}(\vec{a}_i)=\{\vec{x}\in\mathbb{R}^n:|\vec{x}-\vec{a}_i|<\varrho_i\}$ is
the pre-existing spherical hole centered at $\vec{a}_i$ with small radius $\varrho_i>0$,
$W: M^{n \times n}_+\rightarrow\mathbb{R}^+$ is the elastic stored energy density of the
material, $M_+^{n \times n}$ denotes the set of $n \times n $ matrices with positive
determinant, $\Gamma_0$ is the boundary of $\Omega$.

A typical example of the elastic stored energy density is of the form
\begin{equation}
    \label{energy density}
  W(F)=\omega |F|^p+ g(\det F),\quad \forall F \in M^{n\times n}_+,
\end{equation}
where $\omega>0$ is a material constant, $n-1<p<n$,
and $g:(0,\infty)\rightarrow
(0,\infty)$ is a continuously differentiable strictly convex function characterizing
the compressibility of the material and satisfies
\begin{equation}\label{asmp g}
g(d)\rightarrow+\infty\mbox{ as } d \rightarrow 0,\mbox{ and }
\frac{g(d)}{d}\rightarrow+\infty \mbox{ as } d\rightarrow +\infty.
\end{equation}
As was shown by Ball \cite{Ball82}, this kind of functional can have a singular minimizer
displaying cavitation. Further studies on the existence of singular minimizers in
Sobolev spaces are referred to \cite{Henao2009, Muller, Sival2000}.

One of the main difficulties in the computation of immense growth of
voids is the orientation-preservation of the finite element
deformation, which is a crucial constraint and is characterized by the pointwise
positivity of the Jacobian determinant of the deformation gradient.
For the conforming piecewise affine finite element, the condition leads to an
unbearably large amount of degrees of freedom \cite{Xu and
Henao}. Some numerical methods \cite{Lian Dual, Lian and Li
iso, Lian and Li 3, Xu and Henao} have been designed to overcome the difficulty,
and have shown some success numerically. However, strict analytical
results are insufficient. The only practical analytical results for the
cavitation computation known to the authors so far are \cite{SuLiRectan}, where a
sufficient orientation-preservation condition and the interpolation error
estimates were given for a dual-parametric bi-quadratic finite element method,
and \cite{detp}, where a set of sufficient and necessary orientation-preservation
conditions for the quadratic iso-parametric finite element interpolation functions
of radially symmetric cavity deformations are derived.

In this paper, we analyze the approximation properties of
a quadratic iso-parametric finite element for the typical cavitation problem.
The analytical results on the errors of finite element interpolation functions
lead to a delicate relationship between the elastic energy error and the mesh
parameters, which together with the orientation-preservation conditions
(see Remark~\ref{rm-orientation} and \cite{detp})
enable us to establish a mesh distribution strategy guaranteeing that
the corresponding finite element cavitation solution is orientation preserving
and its relative error on the elastic energy is $O(h^2)$, where $h$ is
the mesh size in the far field, {\em i.e.} a given distance away from the cavity.
Above all, for the first time to our knowledge, the convergence of the
finite element cavitation solutions in $W^{1,p}$ norm is proved. In fact, the
numerical experiments show that the optimal order of convergence rate is achieved
by the numerical cavitation solutions obtained on the meshes produced by our
meshing strategy.

Since the cavitation solution is generally considered to be quite regular except
in a neighborhood of the voids, where the material experiences extremely large expansion
dominant deformations and the difficulty of the computation as well as
analysis lies, we restrict ourselves to a simplified problem with
$\Omega_\varrho = B_1(\vec{0}) \setminus B_{\varrho}(\vec{0})$ in $\mathbb{R}^2$
and a simple expansionary boundary condition given by $\vec{u}_0=\lambda \vec{x}$.

The structure of the paper is as follows. In \S~2, we introduce
the iso-parametric finite element method and a radially symmetric large expansion
accommodating triangulation method briefly. \S~3 is
devoted to analyzing the interpolation errors of the cavitation solutions.
The meshing strategy is given in \S~4, where the convergence theorem is also established.
The numerical results are presented
in \S~5. Some concluding remarks are made in \S~6.

\section{Preliminaries}

\subsection{The quadratic iso-parametric FEM}
\begin{figure}[h!]
 \begin{minipage}[l]{0.5\textwidth}
 \centering \includegraphics[width=2in]{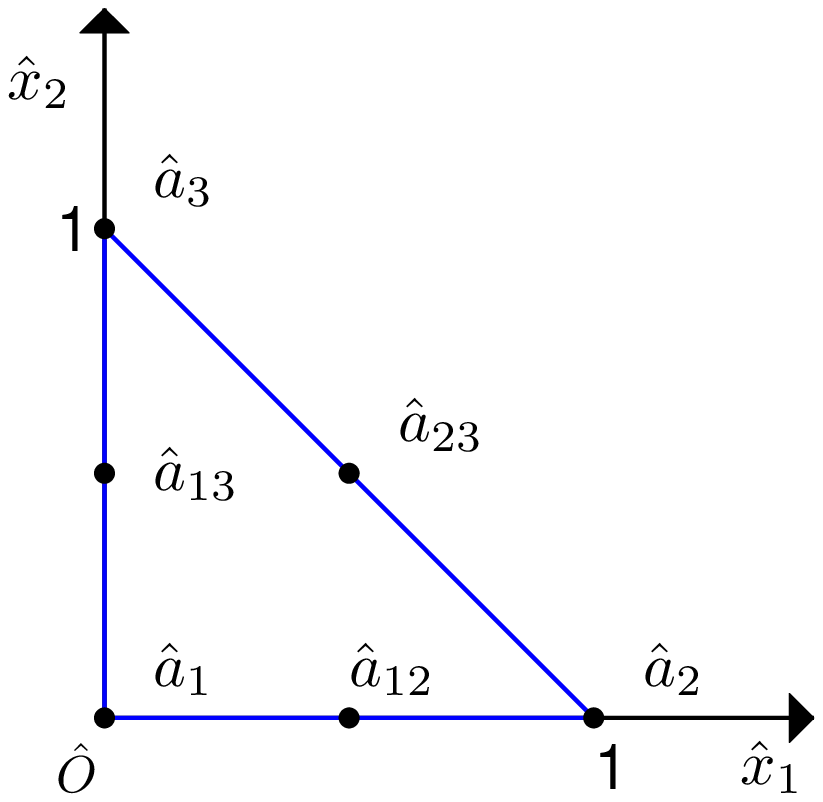}
 \caption{The reference element $\hat{T}$.} \label{dengcan reference element}
\end{minipage}
\begin{minipage}[l]{0.5\textwidth}
 \centering \includegraphics[width=2in]{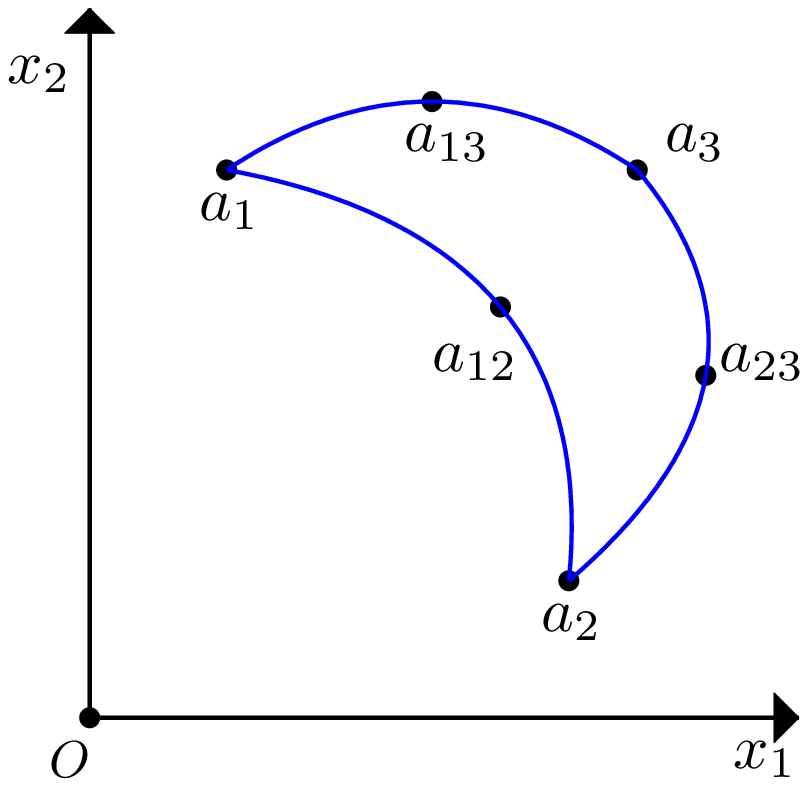}
 \caption{A curved triangular element T.}\label{dengcan curved element}
\end{minipage}
\end{figure}

Let $(\hat{T},\hat{P},\hat{\Sigma})$ be a quadratic Lagrange reference
element. Define $F_T:\hat{T}\rightarrow \mathbb{R}^2$
\begin{equation}\label{e:quadratic_mapping}
        \left\{
                \begin{aligned}
                   & F_T\in \hat{P}^2=(P_2(\hat{T}))^2,\\
                   &\vec{x}=F_T(\hat{\vec{x}})=
                   \sum\limits_{i=1}^{3}\vec{a}_i\hat{\mu}_i(\hat{\vec{x}})+
        \sum\limits_{1\leq i<j\leq 3}\vec{a}_{ij} \hat{\mu}_{ij}(\hat{\vec{x}}),
                \end{aligned}
         \right .
    \end{equation}
where $\vec{a}_i,1\leq i\leq 3$, and $\vec{a}_{ij},1\leq i<j\leq 3$ are given points
in $\mathbb{R}^2$, and
$$
\hat{\mu}_i(\hat{\vec{x}})=\hat{\lambda}_i(\hat{\vec{x}})(2\hat{\lambda}_i
(\hat{\vec{x}})-1), \quad \hat{\mu}_{ij}(\hat{\vec{x}})=
4\hat{\lambda}_i(\hat{\vec{x}})\hat{\lambda}_j(\hat{\vec{x}}),
$$
with $\hat{\lambda}_i(\hat{\vec{x}}),1\leq i\leq 3$ being the barycentric
coordinates of $\hat{T}$.
If the map $F_T$ defined above is an injection, then, $T=F_T(\hat{T})$ is a
curved triangular element as shown in Figure \ref{dengcan curved element}.
The standard quadratic iso-parametric finite element is defined as a finite
element triple $(T,P_T,\Sigma_T)$ with
\begin{equation}
 \label{iso_element}
    \left\{
        \begin{aligned}
  &   T=F_T(\hat{T}) \text{ being a curved triangle element},\\
  &   P_T=\{p:T\rightarrow\mathbb{R}~|~p=\hat{p}\circ F_T^{-1},\;
  \hat{p}\in\hat{P}=P_2(\hat{T})\},\\
  &   \Sigma_T=\{p(\vec{a}_i),1\le i\le 3; p(\vec{a}_{ij}),1\le i<j\le 3\}.
   \end{aligned}
   \right .
\end{equation}

\subsection{Large expansion accommodating triangulations}
\begin{figure}[ht!]
     \begin{minipage}[l]{0.5\textwidth}
 \centering \includegraphics[width=2.35in]{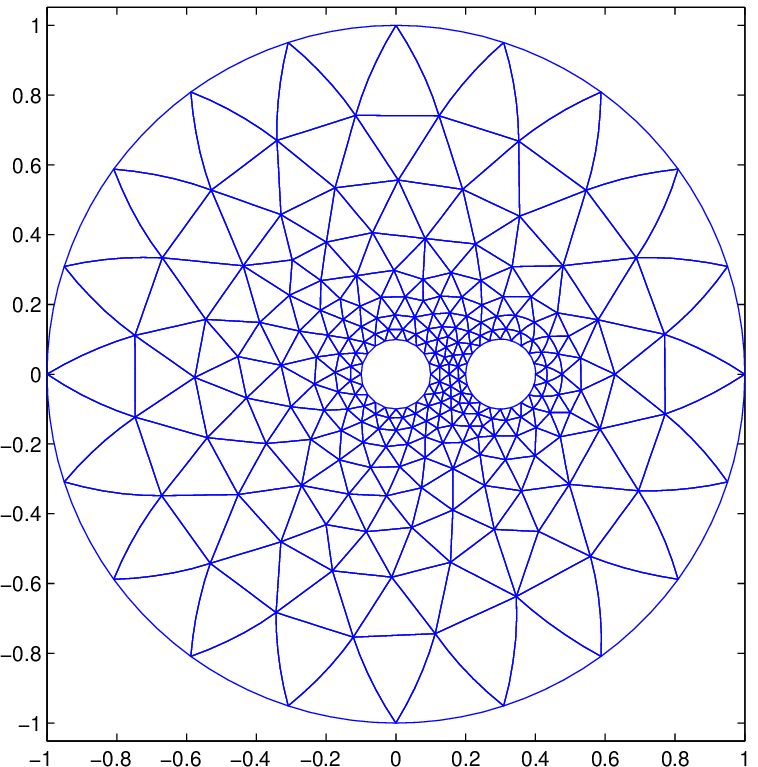}
 \caption{An EasyMesh $\mathcal{J'}$.} \label{easy_mesh}
 \end{minipage}
\begin{minipage}[l]{0.5\textwidth}
 \centering \includegraphics[width=2.35in]{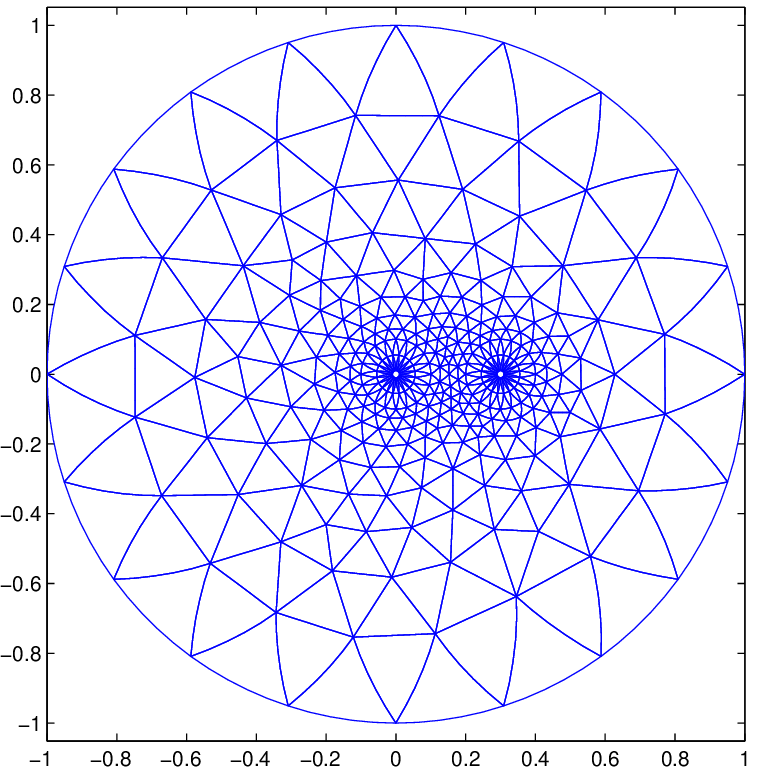}
 \caption{$\mathcal{J}$, a mesh adapted to cavity.} \label{final_mesh}
\end{minipage}
\end{figure}

\begin{figure}[h!]
     \begin{minipage}[c]{0.5\textwidth}\hspace{-35mm}
     \subfigure[A standard layer of curved triangulation.]{
     \includegraphics[width=13.5cm, height=7cm]{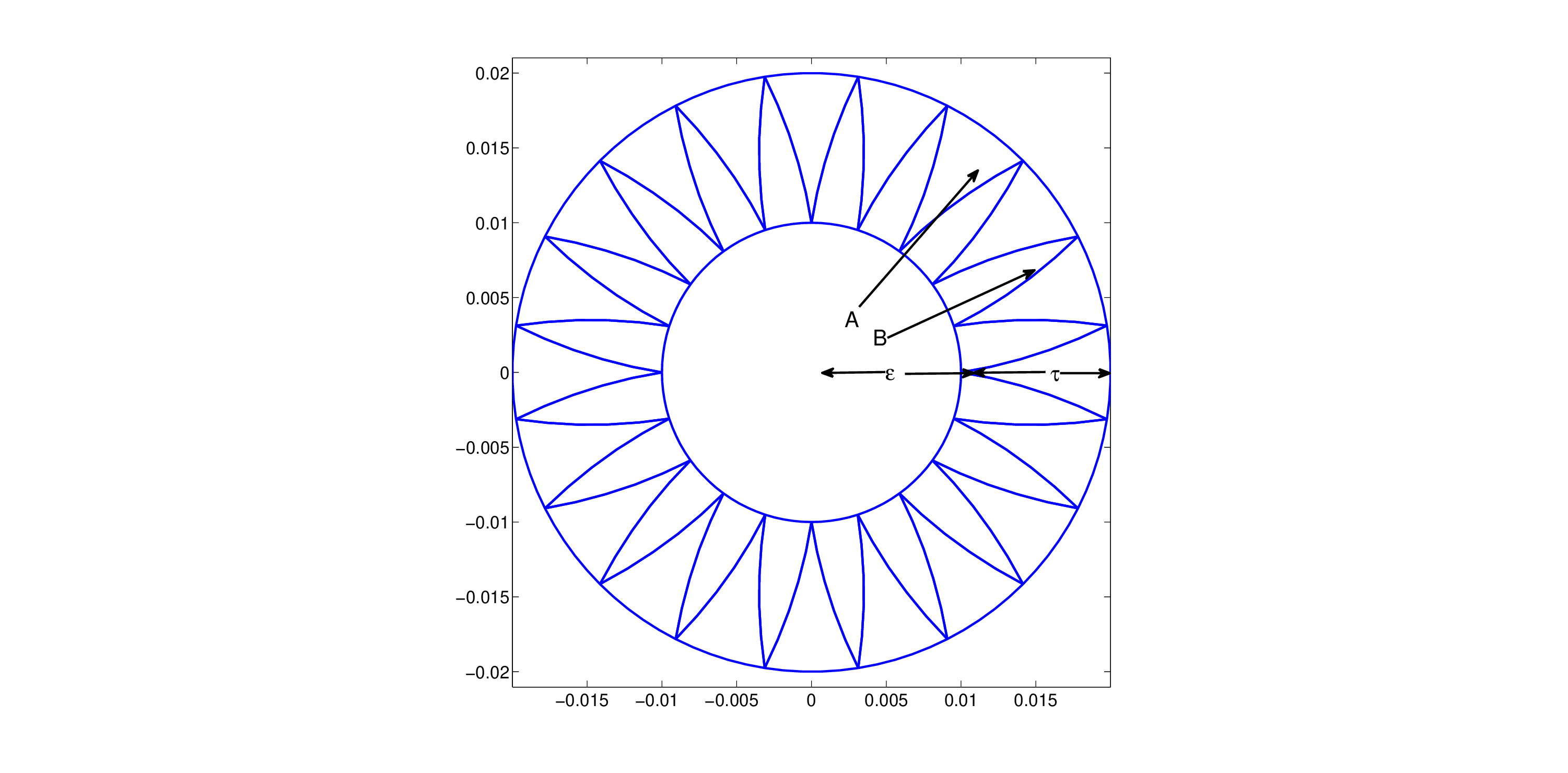}\label{tuli}}
     \end{minipage}
     \begin{minipage}[c]{0.5\textwidth}\hspace{-35mm}
     \subfigure[A conforming layer links 2 standard ones.]{
     \includegraphics[width=13.5cm, height=7cm]{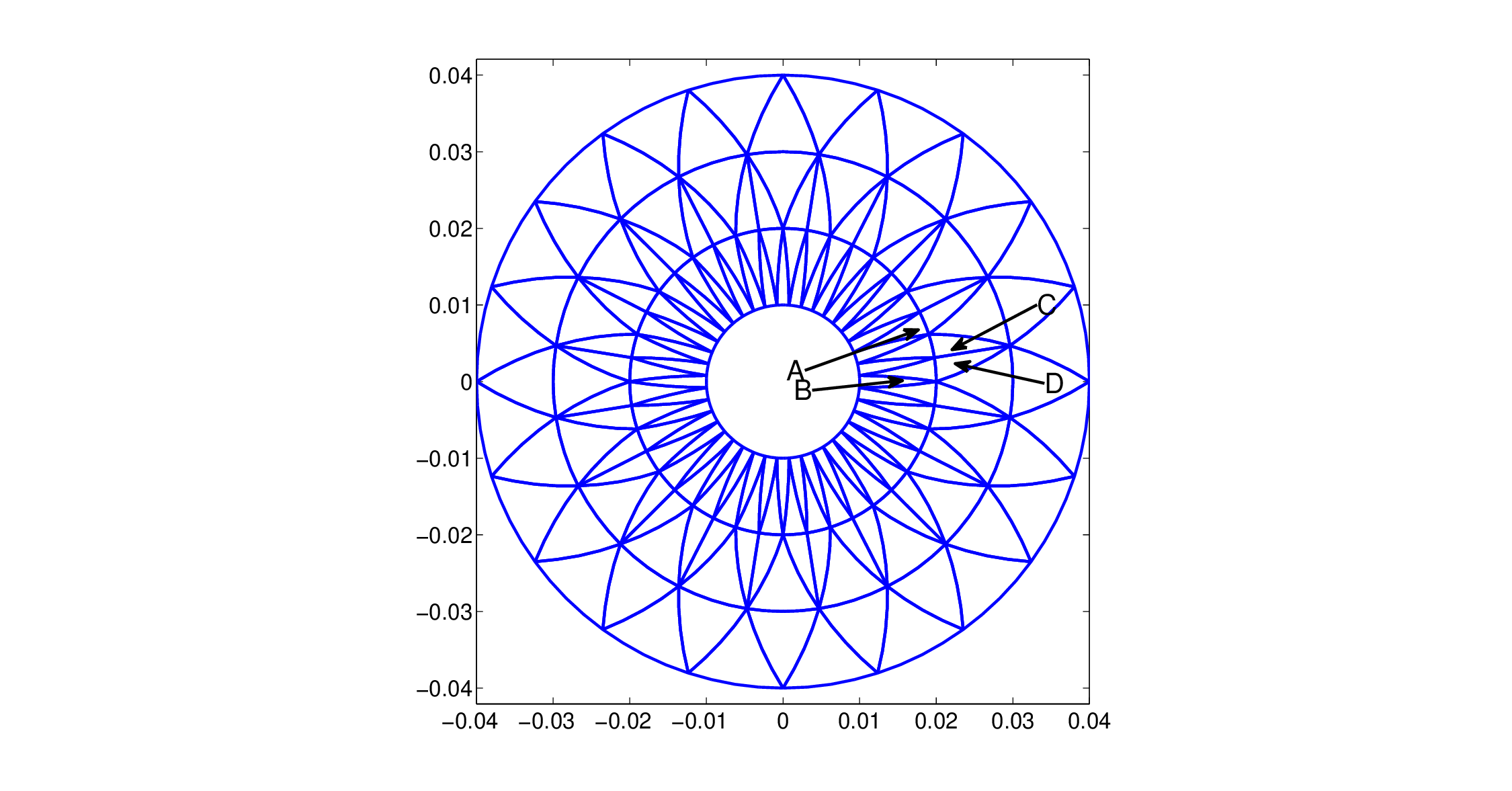}\label{tuli1}}
     \end{minipage}
     \vspace*{-4mm}
 \caption{Triangulation layers with elements of types A, B, C and D.} \label{tuli0}
\end{figure}

Let's have a look at how the iso-parametric FEM is applied in cavitation computation. As
introduced in \cite{Lian and Li iso}, let $\vec{x}_k$ ($k=1,\cdots,m$) be the center of
the $k-$th defect with small radius $\varrho_k$. In the domain far away from the defects
where the deformation is regular, the mesh is given by general straight edged triangulation.
To accommodate the large expansionary deformation
around the defects, the triangulation near the defects is given by the quadratic mapping
as in \eqref{e:quadratic_mapping}, where given 3 vertices $\vec{a}_i$, denote
$(r_k(\vec{x}),\theta_k(\vec{x}))$ the local polar coordinates of $\vec{x}$ with
respect to $\vec{x}_k$, set
\begin{equation}
  \label{eq:make_curved_mesh}
  \vec{a}_{ij}=(r_{ij}\cos\theta_{ij},r_{ij}\sin\theta_{ij})+\vec{x}_k,
\end{equation}
where
\begin{equation*}
 r_{ij}=\frac{r(\vec{a}_i)+r(\vec{a}_j)}{2},\;\;\;
 \theta_{ij}=\frac{\theta(\vec{a}_i)+\theta(\vec{a}_j)}{2}.
 \end{equation*}
Moreover, the elements on the boundary are also adjusted to achieve a better
approximation of the region. With this kind
of curved elements near the defects and general straight triangles elsewhere, the
mesh can better accommodate the locally large expansionary deformations. As an example,
an EasyMesh produced mesh $\mathcal{J'}$ with $m=2$ is shown in
Figure~\ref{easy_mesh}, and the final mesh $\mathcal{J}$ (see Figure~\ref{final_mesh})
is obtained by adding to $\mathcal{J'}$ around each defect two layers of
radially symmetric mesh of the kind shown in Figure~\ref{tuli}, which is a standard curved
triangulation around a prescribed circular ring with inner radius $\epsilon = 0.01$
and thickness $\tau =0.01$. Figure~\ref{tuli1} shows that an outer layer of standard curved
triangulation is connected to a doubly refined inner one by a conforming layer of
nonstandard curved triangulation. For the convenience of reference, we
classify the curved triangular elements in Figure~\ref{tuli0} into four basic types,
and denote them as types A, B, C and D.

\section{Interpolation errors of cavity deformations}
There are standard error estimates \cite{Ciarlet} on the interpolation functions in the
iso-parametric finite element function spaces, however, because of the highly anisotropic
deformation inevitably involved in the cavitation computation, the error bounds so obtained
depend generally on the size of the initial void $\epsilon_0$ and blow up when
$\epsilon_0\rightarrow 0$. In this section, with some more sophisticated manipulations
and calculations, we are able to establish the $\epsilon_0$-independent interpolation
error estimates, including those on the deformation and its
Jacobian determinant, and in particular the elastic energy, for the radially symmetric
cavity deformations in the quadratic iso-parametric finite element function spaces
defined on the meshes consisting of only elements of types A, B, C and D.

To have a better picture in our mind for the problem given below, we first introduce
some notations. Let $\epsilon$ and $\tau$ represent respectively the inner radius
and the thickness of the circular annulus as shown in Figure~\ref{tuli},
let $N$ be the number of the evenly spaced nodes on the outer circle of the circular
ring, and denote
$\Omega_{(\epsilon,\tau)}=\{\vec{x} \in \mathbb{R}^2: \epsilon \le |\vec{x}| \le
\epsilon+\tau\}$. Throughout the paper, the notation
$\Phi \lesssim \Psi$ means that there exists a generic constant $C$ independent of
$\epsilon$ and $\tau$ such that $|\Phi| \le C \Psi$, and $\Phi\sim \Psi$ means that
$\Psi \lesssim \Phi \lesssim \Psi$.

\subsection{The interpolation function and its Jacobian }

For a radially symmetric function
$\vec{v}(\vec{x})=\frac{s(|\vec{x}|)}{|\vec{x}|} \vec{x}$, the iso-parametric
finite element interpolation function can be written as (see \S~2)
\begin{equation}
\Pi \vec{v}(\vec{x})=\sum\limits_{i=1}^{3}\vec{b}_i\hat{\mu}_i(\hat{\vec{x}})+
        \sum\limits_{1\leq i<j\leq 3}\vec{b}_{ij}\hat{\mu}_{ij}(\hat{\vec{x}}),
\label{eq deformation}
\end{equation}
where $\hat{\vec{x}}=F_T^{-1}(\vec{x})$, and where, for
a representative of type A element,
$\vec{b}_1=(s_0,0)$,
$\vec{b}_2=s_1(\cos{\frac{\pi}{N}},-\sin{\frac{\pi}{N}})$,
$\vec{b}_3=s_1(\cos{\frac{\pi}{N}},\sin{\frac{\pi}{N}})$,
$\vec{b}_{12}=s_{1/2}(\cos{\frac{\pi}{2N}},-\sin{\frac{\pi}{2N}})$,
$\vec{b}_{13}=s_{1/2}(\cos{\frac{\pi}{2N}},\sin{\frac{\pi}{2N}})$,
$\vec{b}_{23}=(s_1,0)$.
On this element, let $y=\hat{x}_1+\hat{x}_2$ and $z=\hat{x}_1\hat{x}_2$, we have
\begin{equation}\label{PiA}
\Pi \vec{v}(\vec{x})=\left(s_0+\alpha_2 y+2\alpha_1
y^2-4s_1\sin^2{\frac{\pi}{2N}}(\hat{x}_1^2+\hat{x}_2^2),
(2\gamma y-\beta)(\hat{x}_2-\hat{x}_1)\right),
\end{equation}
where
\begin{eqnarray}
\alpha_1 &=&s_0+s_1-2s_{1/2}\cos{\frac{\pi}{2N}}, \label{pAalpha1} \\
\alpha_2&=&-3s_0-s_1\cos{\frac{\pi}{N}}+4s_{1/2}
\cos{\frac{\pi}{2N}}, \label{pAalpha2}\\
\beta&=&s_1 \sin{\frac{\pi}{N}}-4s_{1/2}
\sin{\frac{\pi}{2N}}, \label{pAbeta} \\
\gamma&=&s_1 \sin{\frac{\pi}{N}}-2s_{1/2}
\sin{\frac{\pi}{2N}}.\label{pAgamma}
\end{eqnarray}
Hence
\begin{equation}\label{JacobipA}
\frac{\partial \Pi \vec{v}}{\partial \hat{\vec{x}}}=\left(
\begin{array}{cc}
\alpha-8s_1\hat{x}_1 \sin^2 {\frac{\pi}{2N}} & \alpha-8s_1\hat{x}_2 \sin^2
{\frac{\pi}{2N}} \\
\beta-4\gamma \hat{x}_1 & -\beta+4 \gamma \hat{x}_2 \\
\end{array}
\right),
\end{equation}
where $\alpha=4\alpha_1 y+\alpha_2$. It follows that
\begin{multline}\label{det PipA}
\det \frac{\partial \Pi \vec{v}}{\partial \hat{\vec{x}}}(\hat{x}_1,\hat{x}_2)=
H(y,z)\triangleq
16 \gamma \alpha_1 y^2-64s_1\gamma \sin^2{\frac{\pi}{2N}}z \\ +
\big(-8\beta \big(\alpha_1-s_1\sin^2{\frac{\pi}{2N}}\big)+
4 \gamma \alpha_2\big)y-2\beta \alpha_2.
\end{multline}

On a representative element of type B with
$\vec{b}_1=(s_1,0)$, $\vec{b}_2=s_0(\cos{\frac{\pi}{N}},\sin{\frac{\pi}{N}})$,
$\vec{b}_3=s_0(\cos{\frac{\pi}{N}},-\sin{\frac{\pi}{N}})$,
$\vec{b}_{12}=s_{1/2}(\cos{\frac{\pi}{2N}},\sin{\frac{\pi}{2N}})$,
$\vec{b}_{13}=s_{1/2}(\cos{\frac{\pi}{2N}},-\sin{\frac{\pi}{2N}})$,
$\vec{b}_{23}=(s_0,0)$, again denote $y=\hat{x}_1+\hat{x}_2$ and
$z=\hat{x}_1\hat{x}_2$, one has
\begin{equation}\label{PiB}
\Pi \vec{v}(\vec{x})=\left(s_1+\bar{\alpha}_2 y+2\bar{\alpha}_1
y^2-4s_0\sin^2{\frac{\pi}{2N}}(\hat{x}_1^2+\hat{x}_2^2), (2\bar{\gamma}
y-\bar{\beta})(\hat{x}_1-\hat{x}_2)\right),
\end{equation}
where
\begin{eqnarray}
\bar{\alpha}_1&=&s_0+s_1-2s_{1/2}\cos{\frac{\pi}{2N}},\label{pBalpha1}\\
\bar{\alpha}_2&=&-3s_1-s_0\cos{\frac{\pi}{N}}+4s_{1/2}
\cos{\frac{\pi}{2N}}, \label{pBalpha2}\\
\bar{\beta}&=&s_0 \sin{\frac{\pi}{N}}-4s_{1/2} \sin{\frac{\pi}{2N}},\label{pBbeta}
\\ \bar{\gamma}&=&s_0 \sin{\frac{\pi}{N}}-2s_{1/2} \sin{\frac{\pi}{2N}}.\label{pBgamma}
\end{eqnarray}
Hence
\begin{equation}\label{JacobipB}
\frac{\partial \Pi \vec{v}}{\partial \hat{\vec{x}}}=\left(
\begin{array}{cc}
\bar{\alpha}-8s_0\hat{x}_1 \sin^2 {\frac{\pi}{2N}} & \bar{\alpha}-8s_0\hat{x}_2 \sin^2
{\frac{\pi}{2N}} \\
-\bar{\beta}+4\bar{\gamma} \hat{x}_1 & \bar{\beta}-4 \bar{\gamma} \hat{x}_2\\
\end{array}
\right),
\end{equation}
\begin{multline}\label{det PipB}
\det \frac{\partial \Pi \vec{v}}{\partial \hat{\vec{x}}}=H(y,z)=
-16 \bar{\gamma} \bar{\alpha}_1 y^2 +64 s_0\bar{\gamma} \sin^2{\frac{\pi}{2N}}z\\
+
(8\bar{\beta} (\bar{\alpha}_1-s_0\sin^2{\frac{\pi}{2N}})-4 \bar{\gamma} \bar{\alpha}_2)y+2
\bar{\beta} \bar{\alpha}_2.
\end{multline}

On a representative element of type C with
$\vec{b}_1=(s_0,0)$, $\vec{b}_2=(s_1,0)$,
$\vec{b}_3=s_0(\cos{\frac{\pi}{N}},\sin{\frac{\pi}{N}})$,
$\vec{b}_{12}=(s_{1/2},0)$,
$\vec{b}_{13}=s_0(\cos{\frac{\pi}{2N}},\sin{\frac{\pi}{2N}})$,
$\vec{b}_{23}=s_{1/2}(\cos{\frac{\pi}{2N}},\sin{\frac{\pi}{2N}})$, one has
\begin{multline}\label{PiC}
\Pi \vec{v}(\vec{x})=\Big(s_0+\tilde{\alpha}_1\hat{x}_1+s_0\tilde{\alpha}_2\hat{x}_2+
2\tilde{\alpha}_3\hat{x}_1^2
-8\sin^2{\frac{\pi}{4N}}\hat{x}_2\big(s_0\cos{\frac{\pi}{2N}}\hat{x}_2-(s_0-s_{1/2})
\hat{x}_1\big), \\
2\sin{\frac{\pi}{2N}}\hat{x}_2\big(s_0(2-\cos{\frac{\pi}{2N}})
-4s_0\sin^2{\frac{\pi}{4N}}\hat{x}_2+2(s_{1/2}-s_0)\hat{x}_1\big) \Big),
\end{multline}
where $\tilde{\alpha}_1=4s_{1/2}-s_1-3s_0$,
$\tilde{\alpha}_2=4\cos{\frac{\pi}{2N}}-\cos{\frac{\pi}{N}}-3$,
$\tilde{\alpha}_3=s_0+s_1-2s_{1/2}$.

Similarly, on a representative element of type D with
$\vec{b}_1=(s_0,0)$, $\vec{b}_2=s_0(\cos{\frac{\pi}{N}},-\sin{\frac{\pi}{N}})$,
$\vec{b}_3=(s_1,0)$,
$\vec{b}_{12}=s_0(\cos{\frac{\pi}{2N}},-\sin{\frac{\pi}{2N}})$,
$\vec{b}_{13}=(s_{1/2},0)$,
$\vec{b}_{23}=s_{1/2}(\cos{\frac{\pi}{2N}},-\sin{\frac{\pi}{2N}})$, one has
\begin{multline}\label{PiD}
\Pi \vec{v}(\vec{x})=\Big(s_0+s_0\tilde{\alpha}_2\hat{x}_1+\tilde{\alpha}_1\hat{x}_2+
2\tilde{\alpha}_3\hat{x}_2^2-8\sin^2{\frac{\pi}{4N}}\hat{x}_1\big(s_0
\cos{\frac{\pi}{2N}}\hat{x}_1-(s_0-s_{1/2})
\hat{x}_2\big),\\
-2\sin{\frac{\pi}{2N}}\hat{x}_1\big(s_0(2-\cos{\frac{\pi}{2N}})
-4s_0\sin^2{\frac{\pi}{4N}}\hat{x}_1+2(s_{1/2}-s_0)\hat{x}_2\big) \Big).
\end{multline}

Throughout this section, we assume that
$\vec{u}(\vec{x})=\frac{r(|\vec{x}|)}{|\vec{x}|}\vec{x}$ represent a cavity deformation,
where $r(R)$ defined on $(0,1]$ is smooth, positive, increasing, and convex with
$r^{(j)}(R)$($j=0, 1, 2, 3$) bounded. Moreover, it satisfies
$\inf\limits_{R\in(0,1]} r(R)>0$, $m R \le r'(R) \le MR$, with $0<m< M$ which
is shown to be naturally satisfied for the energy minimizer of \eqref{energy density}
in \cite{detp}. For simplicity of the notations, we denote $y=\hat{x}_1+\hat{x}_2\in[0,1]$,
$z=\hat{x}_1 \hat{x}_2\in[0,\frac{y^2}{4}]$, where
$\hat{\vec{x}}=(\hat{x}_1,\hat{x}_2)\in \hat{T}$.

\subsection{The error of the interpolation function}
We will estimate, in this subsection, the errors between
$\vec{u}(\vec{x})=\frac{r(|\vec{x}|)}{|\vec{x}|}\vec{x}$ and its interpolation
function $\Pi \vec{u}(\vec{x})$ in the polar coordinates.
Let $(R,\theta)$ be the polar coordinates of $\vec{x}$ on the reference
configuration and let $(\zeta,\varphi)$ be the polar coordinates of the
interpolation function $\Pi \vec{u}(\vec{x})$ of
$\vec{u}(\vec{x})=\frac{r(|\vec{x}|)}{|\vec{x}|}\vec{x}=(r(R),\theta)$ on
the deformed configuration.

\begin{lemma}\label{wgp}
Denote $\vec{x}=F_T(\hat{\vec{x}})=R(\cos \theta, \sin \theta)$, then for
the typical element A, one has
\begin{eqnarray}
R &=&\epsilon+\tau y+O(\tau N^{-2}+\epsilon N^{-4})=(\epsilon+
\tau y)(1+O(N^{-2})),\label{rhoA}\\
\theta &=&\frac{\pi}{N}(\hat{x}_2-\hat{x}_1)+O(N^{-3}),\label{thetaA}\\
\det \nabla \vec{x}&=&4\tau \sin \frac{\pi}{2N}(\epsilon+\tau y+2\tau
\sin^2 \frac{\pi}{4N})\big(1+O(N^{-2}+\frac{\epsilon}{\tau}N^{-4})\big); \label{detA}
\end{eqnarray}
for the typical element B, we have
\begin{eqnarray}
&&R =\epsilon+\tau (1-y)+O(\tau +\epsilon N^{-2})N^{-2}=(\epsilon+\tau (1-y))(1+O(N^{-2})),
\qquad \label{rhoB}\\
&&\theta =\frac{\pi}{N}(\hat{x}_1-\hat{x}_2)+O(N^{-3}),\;\;\label{thetaB}\\
&&\det \nabla \vec{x}=4\tau \sin \frac{\pi}{2N}\big(\epsilon+\tau(1-y)\big)
\big(1+O(N^{-2}+\frac{\epsilon}{\tau}N^{-4})\big); \label{detB}
\end{eqnarray}
for the typical element C, one has
\begin{eqnarray}
R &=&\epsilon+\tau \hat{x}_1+O(\tau N^{-2}+\epsilon N^{-4})=(\epsilon+
\tau \hat{x}_1)(1+O(N^{-2})),\label{rhoC}\\
\theta &=&\frac{\pi}{N}\hat{x}_2+O(N^{-3}),\label{thetaC}\\
\det \nabla \vec{x}&=&2\tau \sin \frac{\pi}{2N}(\epsilon+\tau \hat{x}_1)
\big(1+O(N^{-2})\big); \label{detC}
\end{eqnarray}
while for the typical element D, one has
\begin{eqnarray}
R &=&\epsilon+\tau \hat{x}_2+O(\tau N^{-2}+\epsilon N^{-4})=(\epsilon+
\tau \hat{x}_2)(1+O(N^{-2})),\label{rhoD}\\
\theta &=&-\frac{\pi}{N}\hat{x}_1+O(N^{-3}),\label{thetaD}\\
\det \nabla \vec{x}&=&2\tau \sin \frac{\pi}{2N}(\epsilon+\tau \hat{x}_2)
\big(1+O(N^{-2})\big). \label{detD}
\end{eqnarray}
\end{lemma}
\textbf{Proof.} For a typical type A element, it follows from \eqref{PiA} that
$x_1=t_1+t_2z, x_2=t_3(\hat{x}_2-\hat{x}_1)$, where
\begin{eqnarray}
t_1&=&\epsilon+\hat{\alpha}_2y+\hat{\alpha}_4 y^2, \qquad\;\;
t_2=2\hat{\gamma}y-\hat{\beta}, \quad\;
t_3=8(\epsilon+\tau)\sin^2{\frac{\pi}{2N}},\nonumber\\
\hat{\alpha}_4&=&2\epsilon+2(\epsilon+\tau)\cos{\frac{\pi}{N}}-
4(\epsilon+\frac{\tau}{2})\cos{\frac{\pi}{2N}}\nonumber\\
&=&-2\epsilon \sin^2\frac{\pi}{2N}+O(\tau N^{-2}+\epsilon N^{-4}),\label{halpha4}
\end{eqnarray}
$\hat{\alpha}_2$, $\hat{\beta}$, $\hat{\gamma}$ are given by
\eqref{pAalpha2}-\eqref{pAgamma} by taking $s(t)=t$, i.e.,
\begin{eqnarray}
\hat{\alpha}_2&=&-3\epsilon+4(\epsilon+\frac{\tau}{2})\cos{\frac{\pi}{2N}}-
(\epsilon+\tau) \cos{\frac{\pi}{N}}=\tau+O(\tau N^{-2}+\epsilon N^{-4}),
\qquad \label{halpha2} \\
\hat{\beta}&=&(\epsilon+\tau) \sin{\frac{\pi}{N}}-
4(\epsilon+\frac{\tau}{2}) \sin{\frac{\pi}{2N}}\nonumber \\
&=&-2\epsilon \sin{\frac{\pi}{2N}}-4(\epsilon+\tau)\sin{\frac{\pi}{2N}}
\sin^2{\frac{\pi}{4N}}, \label{hbeta} \\
\hat{\gamma}&=&(\epsilon+\tau) \sin{\frac{\pi}{N}}-
2(\epsilon+\frac{\tau}{2}) \sin{\frac{\pi}{2N}}\nonumber \\
&=&\tau\sin{\frac{\pi}{2N}}-4(\epsilon+\tau)\sin{\frac{\pi}{2N}}
\sin^2{\frac{\pi}{4N}}. \label{hgamma}
\end{eqnarray}
Thus, by $0\le y \le 1$, $\frac{y}{\epsilon+\tau y}\le \frac{1}{\epsilon+\tau}$,
one has
\begin{eqnarray*}
x_1&=&\epsilon+\tau y-2\epsilon \sin^2 \frac{\pi}{2N}(y^2-4z)+O(\tau +
\epsilon N^{-2})N^{-2}y\\
&=&(\epsilon+\tau y)(1+O(N^{-2})),\\
x_2&=&2\sin \frac{\pi}{2N}(\epsilon+\tau y-2(\epsilon+\tau)
\sin^2 \frac{\pi}{4N}(2y-1))(\hat{x}_2-\hat{x}_1),
\end{eqnarray*}
hence
\begin{eqnarray*}
R^2&=&(\epsilon+\tau y)^2-4\epsilon(\epsilon+\tau y)
\sin^2 \frac{\pi}{2N}(y^2-4z)+(\epsilon+\tau y)O(\epsilon N^{-4}+\tau N^{-2})y\\
&&\qquad \quad +\tau^2 N^{-6}y^2+4(\epsilon+\tau y)^2 \sin^2 \frac{\pi}{2N}(y^2-4z)\\
&=&(\epsilon+\tau y)^2\Big(1+\frac{O(\tau N^{-2}+\epsilon N^{-4})y}{\epsilon+\tau y}
+\frac{O(\tau^2 N^{-6})y^2}{(\epsilon+\tau y)^2}\Big).
\end{eqnarray*}
This gives \eqref{rhoA}. On the other hand,
\begin{eqnarray*}
|\tan \theta-\frac{\pi}{N}(\hat{x}_2-\hat{x}_1)|
&=&|\frac{x_2-\frac{\pi}{N}x_1(\hat{x}_2-\hat{x}_1)}{x_1}|\\
&\lesssim & \frac{(\epsilon+\tau y)(2\sin \frac{\pi}{2N}-\frac{\pi}{N})
+O((\epsilon+\tau)N^{-3})}{(\epsilon+\tau y)(1+O(N^{-2}))}y\\
&\lesssim&N^{-3},
\end{eqnarray*}
which gives \eqref{thetaA}. Next consider $\det \nabla \vec{x}$. It follows from
\eqref{pAalpha2}-\eqref{pAgamma}, \eqref{det PipA} and $z\le y^2/4 \le y/4$ that
$$
\det \nabla \vec{x}=2t_2 \hat{\alpha}_2+O((\epsilon+\tau)\tau N^{-3}+\epsilon^2N^{-5})y.
$$
Note that $\hat{\alpha}_2=\tau(1+O(N^{-2}+\frac{\epsilon}{\tau}N^{-4}))$,
$t_2=2(\epsilon+\tau y+2\tau \sin^2 \frac{\pi}{4N})\sin \frac{\pi}{2N}(1+O(N^{-2}))$,
thus \eqref{detA} follows.

For the type B element, by \eqref{pBalpha1}-\eqref{pBgamma},
\begin{eqnarray*}
x_1&=&\epsilon+(\tau-4\tau \sin^2\frac{\pi}{4N}y-8\epsilon \sin^4 \frac{\pi}{4N}y)(1-y)-
2\epsilon \sin^2\frac{\pi}{2N}(y^2-4z),\\
x_2&=&\big(-2\sin \frac{\pi}{2N}(\epsilon+\tau (1-y))-4\epsilon \sin \frac{\pi}{2N}
\sin^2 \frac{\pi}{4N}(1-2y)\big)(\hat{x}_2-\hat{x}_1),
\end{eqnarray*}
which gives \eqref{rhoB}, \eqref{thetaB}. While by \eqref{det PipB},
\begin{eqnarray*}
&&\det \nabla \vec{x}\\
&=&4\tau \sin \frac{\pi}{2N}\big((\epsilon+\tau (1-y))+(2y^2-3y+1)O(\tau N^{-2})+
O(\epsilon N^{-2}+\frac{\epsilon^2}{\tau}N^{-4})\big)\\
&=&4\tau \sin \frac{\pi}{2N}(\epsilon+\tau (1-y))\big(1+O(N^{-2}+
\frac{\epsilon}{\tau}N^{-4})+\frac{2y^2-3y+1}{\epsilon+\tau(1-y)}
O(\tau N^{-2})\big).
\end{eqnarray*}
Since $\sup\limits_{0\le y\le 1}|\frac{2y^2-3y+1}{\epsilon+\tau(1-y)}|=
\sup\limits_{0\le y\le 1}|\frac{(2y-1)y}{\epsilon+\tau y}|\le \frac{1}{\epsilon+\tau}$,
then we obtain \eqref{detB}.

For the type C element, by \eqref{PiC},
\begin{eqnarray*}
x_1&=&\epsilon+\tau \hat{x}_1-8\epsilon \sin^4 {\frac{\pi}{4N}}\hat{x}_2
-8\epsilon \cos {\frac{\pi}{2N}}\sin^2 {\frac{\pi}{4N}}\hat{x}_2^2
-4\tau \sin^2 \frac{\pi}{4N}\hat{x}_1\hat{x}_2,\\
x_2&=&2\sin \frac{\pi}{2N}\hat{x}_2(\epsilon+\tau \hat{x}_1)+
4\epsilon \sin \frac{\pi}{2N}\sin^2 \frac{\pi}{4N}\hat{x}_2(1-2\hat{x}_2),
\end{eqnarray*}
which gives \eqref{rhoC}, \eqref{thetaC}. Meanwhile, one has
$$\det \nabla \vec{x}=2\tau \sin \frac{\pi}{2N}\big(\epsilon+\tau \hat{x}_1
+2\epsilon \sin^2 {\frac{\pi}{4N}}(1-6\hat{x}_2+8\hat{x}_2^2)\big),$$
which leads to \eqref{detC}.

Similar arguments yield \eqref{rhoD}-\eqref{detD} for the type D element.
\hfill $\square$

\begin{theorem}\label{error of interpolation}
The error between a cavity deformation
$\vec{u}(\vec{x})=\frac{r(|\vec{x}|)}{|\vec{x}|}\vec{x}=(r(R),\theta)$ and its
interpolation function $\Pi \vec{u}$ satisfies
\begin{eqnarray}\label{error of theta}
\theta - \varphi &=& O(\tau^2N^{-1}+N^{-3}), \\
\label{error of r}
r(R)-\zeta &=& O(\tau^3+\epsilon\tau N^{-2}+N^{-4}),
\end{eqnarray}
where the constants in $O(\cdot)$ depend on $\|r\|_\infty$, $\|r''\|_\infty$,
$\|r^{(3)}\|_\infty$, $\inf r(R)$, $\sup \frac{r'(R)}{R}$.
\end{theorem}
\textbf{Proof.} For a typical type A element as used above, denote
$\vec{X}=\Pi \vec{u}(\vec{x})=(X_1,X_2)$, where $\vec{x}=F_T(\hat{\vec{x}})$.
Some tedious manipulation yields that
$X_1=T_1+T_2z, X_2=T_3(\hat{x}_2-\hat{x}_1)$,
where
\begin{eqnarray*}
&&T_1=r(\epsilon)+\alpha_2 y+\alpha_4 y^2, \;\;\; T_2=2\gamma y-\beta, \;\;\;
T_3=8r(\epsilon+\tau)\sin^2{\frac{\pi}{2N}},\\
&&\alpha_4=-4r(\epsilon+\frac{\tau}{2})\cos{\frac{\pi}{2N}}+2r(\epsilon)+
2r(\epsilon+\tau)\cos{\frac{\pi}{N}},\\
&&\hat{\alpha}_4=2\epsilon+2(\epsilon+\tau)\cos{\frac{\pi}{N}}-
4(\epsilon+\frac{\tau}{2})\cos{\frac{\pi}{2N}},
\end{eqnarray*}
and where $\alpha_2$, $\beta$, $\gamma$ are given by \eqref{pAalpha2}-\eqref{pAgamma}.
By the Taylor expansion,
\begin{eqnarray}
\alpha_2&=&r'(\epsilon)\tau+O(\tau^3+\tau^2N^{-2}+\epsilon\tau N^{-2}+N^{-4}),\label{alpha2}\\
\alpha_4&=&\frac{r''(\epsilon)}{2}\tau^2-2r(\epsilon)\sin^2{\frac{\pi}{2N}}+
O(\tau^3+\tau^2N^{-2}+\epsilon \tau N^{-2}+N^{-4}),\label{alpha4}\\
\beta&=& -2r(\epsilon)\sin{\frac{\pi}{2N}}+O(\tau^2 N^{-1}+N^{-3}),\label{beta}\\
\gamma&=&r'(\epsilon)\tau \sin{\frac{\pi}{2N}}+O(\tau^2 N^{-1}+N^{-3}).\label{gamma}
\end{eqnarray}
Hence
\begin{eqnarray*}
 T_1&=&r(\epsilon+\tau y)-2r(\epsilon)y^2\sin^2{\frac{\pi}{2N}}+
O(\tau^{3}+\tau^2N^{-2}+\epsilon \tau N^{-2}+N^{-4}),\\
T_2&=&2r(\epsilon+\tau y)\sin{\frac{\pi}{2N}}+O(\tau^2N^{-1}+N^{-3}),\\
\tan{\varphi} &=&\frac{\hat{x}_2-\hat{x}_1}{T_1+T_3z}T_2\\
&=&\frac{r(\epsilon+\tau y)\frac{\pi}{N}(1+O(\tau^2+N^{-2}))}{r(\epsilon+
\tau y)(1+O(\tau^3+N^{-2}))}(\hat{x}_2-\hat{x}_1)\\
&=&\frac{\pi}{N}(\hat{x}_2-\hat{x}_1)+O(\tau^2N^{-1}+N^{-3}).
\end{eqnarray*}
This together with \eqref{thetaA} gives \eqref{error of theta}. Thus
\eqref{error of theta} is established on the type A element.

Next we estimate $r(R) - \zeta$. Expand
$\zeta^2=X_1^2+X_2^2=(T_1+T_3z)^2+T_2^2(\hat{x}_2-\hat{x}_1)^2$ as follows
\begin{eqnarray*}
\zeta^2&=&(r(\epsilon+\tau y)-2r(\epsilon)(y^2-4z)\sin^2{\frac{\pi}{2N}}+
O(\epsilon \tau N^{-2}+\tau^3+\tau^2N^{-2}+N^{-4}))^2\\
& &\ \ \ + (2r(\epsilon+\tau y)\sin{\frac{\pi}{2N}}+O(\tau^2N^{-1}+N^{-3}))^2(y^2-4z)\\
&=&r^2(\epsilon+\tau y)+4r(\epsilon+\tau y)(r(\epsilon+\tau y)-r(\epsilon))
(y^2-4z)\sin^2{\frac{\pi}{2N}}+\\
& & \ \ \ O(\epsilon \tau N^{-2}+\tau^3+\tau^2N^{-2}+N^{-4})\\
&=&r^2(\epsilon+\tau y)+O(\epsilon \tau N^{-2}+\tau^3+\tau^2N^{-2}+N^{-4}).
\end{eqnarray*}
Hence it follows from \eqref{rhoA} and $r'(R)\le MR$ that \eqref{error of r} holds on
the type A element.

Similarly, we can show that \eqref{error of theta}-\eqref{error of r} hold on elements
of types B, C and D. \hfill $\square$
\vspace{3mm}

\subsection{The error on the Jacobian determinant}

\begin{theorem}\label{error det}
The error between the Jacobian
determinants of a cavity deformation
$\vec{u}(\vec{x})=\frac{r(|\vec{x}|)}{|\vec{x}|}\vec{x}=(r(R),\theta)$ and
its interpolation function $\Pi \vec{u}$ satisfies
\begin{eqnarray}
\det\frac{\partial \Pi \vec{u}}{\partial \vec{x}}(\hat{\vec{x}})-
\det \frac{\partial \vec{u}}{\partial \vec{x}}(\hat{\vec{x}})&=&
\frac{O(\tau^3N^{-1}+(\epsilon+\tau)\tau N^{-3} +N^{-5})}{\det \nabla \vec{x}}
\label{det hPiu}\\
&=&O\Big(N^{-2}+\frac{\tau^2}{\epsilon}+\frac{N^{-4}}{\epsilon \tau}\Big).\label{eq error det}
\end{eqnarray}
\end{theorem}
\textbf{Proof.} For the representative element of type A,
it follows from \eqref{det PipA}, \eqref{alpha2}-\eqref{gamma} and
$$
\alpha_1-r(\epsilon+\tau)\sin^2{\frac{\pi}{2N}}=
r''(\epsilon+\frac{\tau}{2})\frac{\tau^2}
{4}+O(\tau^4+(\epsilon+\tau)\tau N^{-2}+N^{-4}),
$$
that
\begin{eqnarray*}
& & \hspace*{-8mm}\det \frac{\partial \Pi \vec{u}}{\partial \hat{\vec{x}}}=
-2\beta(\alpha_2+4(\alpha_1- r(\epsilon+\tau)\sin^2{\frac{\pi}{2N}})y)
+4\gamma \alpha_2 y \\ && \ \ \ \ \ \ \ \ \ \ + \ O(\tau^4 N^{-1}+
\epsilon \tau^3 N^{-1}+\tau^2N^{-3}+\epsilon \tau N^{-3}+N^{-5}) \\
&&=4\tau \sin{\frac{\pi}{2N}}r(\epsilon+\tau y)r'(\epsilon+\tau y)+
O(\tau^3N^{-1}+(\epsilon+\tau)\tau N^{-3} +N^{-5}),
\end{eqnarray*}
which gives \eqref{det hPiu}.
Hence by \eqref{detA},
\vspace*{-3mm}
\begin{eqnarray*}\hspace*{-8mm}
\det \frac{\partial \Pi \vec{u}}{\partial \vec{x}}(\hat{\vec{x}})&=&
\frac{r(\epsilon+\tau y)r'(\epsilon+\tau y)+O(\tau^2+(\epsilon+\tau) N^{-2}+
\tau^{-1}N^{-4})} {(\epsilon+\tau y)\big(1+O(\frac{\epsilon+\tau}{\epsilon}
N^{-2}+\frac{\epsilon}{\tau}N^{-4})\big)} \\
&=&\frac{r(\epsilon+\tau y)r'(\epsilon+\tau y)}{\epsilon+\tau y}+
O\Big(\frac{\tau^2}{\epsilon}+\frac{\epsilon+\tau}{\epsilon}N^{-2}
+\frac{N^{-4}}{\epsilon \tau}\Big).
\end{eqnarray*}
On the other hand , it follows from \eqref{rhoA} that
$$
\det \frac{\partial \vec{u}}{\partial \vec{x}}(\hat{\vec{x}})=
\frac{r(R)r'(R)}{R}=
\frac{r(\epsilon+\tau y)r'(\epsilon+\tau y)}{\epsilon+\tau y}+
O(N^{-2}).
$$
Hence
$$
\det\frac{\partial \Pi \vec{u}}{\partial \vec{x}}(\hat{\vec{x}})-
\det \frac{\partial \vec{u}}{\partial \vec{x}}(\hat{\vec{x}})=
O\Big(\frac{\tau^2}{\epsilon}+\frac{\epsilon+\tau}{\epsilon}N^{-2}+
\frac{N^{-4}}{\epsilon \tau}\Big).
$$
Since $\frac{\tau^2}{\epsilon}+\frac{N^{-4}}{\epsilon \tau}\ge
2\sqrt{\frac{\tau}{\epsilon^2}N^{-4}}=2\frac{\sqrt{\tau}}{\epsilon}N^{-2}>
\frac{\tau}{\epsilon}N^{-2}$,
thus \eqref{eq error det} holds on type A elements.

The proof for types B, C, D elements is similar.
\hfill $\square$

\begin{remark}\label{rm-orientation}
The results above imply that, if there exits a constant $d>0$ such that
$\det \frac{\partial \vec{u}}{\partial \vec{x}}\ge d$ which is in fact the case for the
radially symmetric cavity solution, then there exist $C_1>0$,
$C_2>0$, $C_3>0$ such that $\det \frac{\partial \Pi \vec{u}}{\partial \vec{x}}>0$
as long as $\tau\le C_1\epsilon^{1/2}$ and $N^{-1}\le \min\{C_2(\epsilon \tau)^{1/4}, C_3\}$,
which is in a good agreement with the result obtained in \cite{detp} for the
orientation-preservation condition.
\end{remark}

\subsection{The error on the elastic energy}

Let $\mathcal{J}(\Omega_{(\epsilon,\tau)})$ be a quadratic iso-parametric
finite element triangulation of $\Omega_{(\epsilon,\tau)}$ consisting of either a
layer of evenly distributed elements of types A and B as shown in Figure~\ref{tuli},
or a layer of evenly distributed elements of types A, C and D as shown in
Figure~\ref{tuli1}, and the corresponding elements are denoted by $T_A$, $T_B$,
$T_C$ and $T_D$ accordingly. Let $\Omega^{\mathcal{J}}_{(\epsilon,\tau)} = \cup_{T \in
\mathcal{J}(\Omega_{(\epsilon,\tau)})} T$. For the energy density function of the
form \eqref{energy density}, $S\subset \mathbb{R}^2$, denote
\begin{eqnarray}
E_1(\vec{u};S) &=& \omega \int_S \left| \nabla \vec{u}\right|^p d \vec{x}, \label{E1}\\
E_2(\vec{u};S) &=& \int_S g\left(\det \nabla \vec{u} \right) d \vec{x}, \label{E2}\\
A(\epsilon,\tau)&=&(2-p)\int_\epsilon^{\epsilon+\tau}t^{1-p}
dt=(\epsilon+\tau)^{2-p}-\epsilon^{2-p}.\label{A}
\end{eqnarray}
Let $E_i(\Pi \vec{u};\Omega^{\mathcal{J}}_{(\epsilon,\tau)})$,
$E(\Pi \vec{u};\Omega^{\mathcal{J}}_{(\epsilon,\tau)})$ be the corresponding counterparts
of the elastic energy of $\Pi \vec{u}$ on the elements in $\Omega_{(\epsilon,\tau)}$.
Then, we have the following result.

\begin{theorem}\label{error energy}
The elastic energies
of a cavity deformation $\vec{u}(\vec{x})=\frac{r(|\vec{x}|)}{|\vec{x}|}\vec{x}$ and its
interpolation function $\Pi \vec{u}$ satisfy
\begin{eqnarray}
&&E_1(\Pi \vec{u};\Omega^{\mathcal{J}}_{(\epsilon,\tau)})=
E_1(\vec{u};\Omega^{\mathcal{J}}_{(\epsilon,\tau)})(1+
O(\tau^2+N^{-2}+\frac{\epsilon}{\tau}N^{-4})),
\label{error E13}\\
&&E_2(\Pi \vec{u};\Omega^{\mathcal{J}}_{(\epsilon,\tau)})=
E_2(\vec{u};\Omega^{\mathcal{J}}_{(\epsilon,\tau)})+O(\tau^3+(\epsilon+\tau)\tau
N^{-2}+N^{-4}), \label{error E23}\\
&&E(\Pi \vec{u};\Omega^{\mathcal{J}}_{(\epsilon,\tau)})=
E(\vec{u};\Omega^{\mathcal{J}}_{(\epsilon,\tau)})(1+O(\tau^2+N^{-2}+\frac{N^{-4}}
{\tau})), \label{error E3}\\
&&E(\vec{u};\Omega^{\mathcal{J}}_{(\epsilon,\tau)})=
E(\vec{u};\Omega_{(\epsilon,\tau)})(1+O(\tau^2+N^{-2}+\frac{\epsilon}{\tau}N^{-4}))
\sim A(\epsilon,\tau). \label{Eo}
\end{eqnarray}
\end{theorem}
\textbf{Proof.} On an element of type A, it follows from \eqref{JacobipA} that
\begin{equation*}
\nabla \Pi \vec{u}(\vec{x})=\frac{\partial \Pi \vec{u}}
{\partial \hat{\vec{x}}} (\nabla \vec{x})^{-1} = \frac{1}{\det \nabla \vec{x}}
\begin{pmatrix} A_{11} & A_{12}\\ A_{21} & A_{22} \end{pmatrix},
\end{equation*}
where by \eqref{alpha2}-\eqref{gamma}, one has
\begin{eqnarray*}
A_{11}&=&4\tau\sin{\frac{\pi}{2N}}(r'(\epsilon)(\epsilon+\tau y)+r''(\epsilon+\tau/2)\epsilon
\tau y)+O(\tau^3N^{-1}+\tau N^{-3}+\epsilon N^{-5})\\
&=& 4r'(\epsilon+\tau y)(\epsilon+\tau y)\tau \sin{\frac{\pi}{2N}}+
O(\tau^3N^{-1}+\tau N^{-3}+\epsilon N^{-5}),\\
A_{12}&=&O(\tau N^{-2}), \\ A_{21}&=&O(\tau N^{-2}),\\
A_{22}&=&4r(\epsilon+\tau y)\tau \sin{\frac{\pi}{2N}}+O(\tau^3N^{-1}+\tau N^{-3}+
\epsilon N^{-5}),
\end{eqnarray*}
where $y=\hat{x}_1+\hat{x}_2$, $\vec{\hat{x}}=F_T^{-1}(\vec{x})$.
Thus, denote $\xi_A(y)=\epsilon+\tau y$ and
$\Upsilon(y)=r^2(\xi_A(y))+r'(\xi_A(y))^2\xi_A(y)^2$, one has
$$
\Big|\frac{\partial \Pi \vec{u}}{\partial \vec{x}}\Big|^2=
\frac{1}{(\det \nabla \vec{x})^2}\sum\limits_{i,j}A_{ij}^2
=\frac{(4\tau\sin \frac{\pi}{2N})^2}{(\det \nabla \vec{x})^2}\Upsilon(y)(1+\iota_1),
$$
where $\iota_1=O(\tau^2+N^{-2}+\frac{\epsilon}{\tau}N^{-4})$. Hence
$$
\Big|\frac{\partial \Pi \vec{u}}{\partial \vec{x}}\Big|^p=
\frac{(4\tau
\sin{\frac{\pi}{2N}})^p}{(\det \nabla \vec{x})^p}\Upsilon(y)^{p/2}(1+\iota_1).
$$
It follows from Lemma~\ref{wgp} or more precisely \eqref{detA} that
\begin{eqnarray*}
&&\int_{T_A}\Big|\frac{\partial \Pi \vec{u}}{\partial \vec{x}}\Big|^pd\vec{x}=\int_{\hat{T}}
\Big|\frac{\partial \Pi \vec{u}}{\partial \vec{x}}\Big|^p\det \nabla \vec{x}d\hat{\vec{x}} \\
&=&4\tau\sin{\frac{\pi}{2N}}\int_{\hat{T}} \Upsilon(y)^{p/2}(\xi_A(y)+
2\tau \sin^2\frac{\pi}{4N})^{1-p}
d\hat{\vec{x}}(1+\iota_1)\\
&=&4\tau\sin{\frac{\pi}{2N}}\int_0^1 y\Upsilon(y)^{p/2}(\xi_A(y)+
2\tau \sin^2\frac{\pi}{4N})^{1-p}
dy(1+\iota_1),
\end{eqnarray*}
where $\hat{T}$ is the reference element as Figure \ref{dengcan reference element}.
Similarly, by Lemma~\ref{wgp} or more precisely \eqref{detB}, on the element of type B,
one has
\begin{eqnarray}
\int_{T_B}\Big|\frac{\partial \Pi \vec{u}}{\partial
\vec{x}}\Big|^pd\vec{x}&=&4\tau\sin{\frac{\pi}{2N}}\int_0^1
y\Upsilon(1-y)^{p/2}\xi_A(1-y)^{1-p}
dy(1+\iota_1)\nonumber\\
&=&4\tau\sin{\frac{\pi}{2N}}\int_0^1 (1-y)\Upsilon(y)^{p/2}\xi_A(y)^{1-p}
dy(1+\iota_1).\label{TBenergy}
\end{eqnarray}
Thus
$$
\int_{T_A\cup T_B}\left| \frac{\partial \Pi \vec{u}}{\partial \vec{x}} \right|^p
d\vec{x}=4\tau \sin{\frac{\pi}{2N}}\int_0^1
\Upsilon(y)^{p/2}\xi_A(y)^{1-p}dy(1+\iota_1) +\digamma,
$$
where by the median formula,
\begin{eqnarray*}
|\digamma|&=&
4\tau \sin{\frac{\pi}{2N}}\int_0^1y\Upsilon(y)^{p/2}\Big(\xi_A(y)^{1-p}-(\xi_A(y)+
2\tau \sin^2\frac{\pi}{4N})^{1-p}\Big)dy\\
&\le &4(p-1)\tau \sin{\frac{\pi}{2N}}\int_0^1y \Upsilon(y)^{p/2}\xi_A(y)^{-p}
2\tau \sin^2\frac{\pi}{4N} dy \\
&< &4(p-1)\tau \sin{\frac{\pi}{2N}}\sin^2\frac{\pi}{4N}\int_0^1
\Upsilon(y)^{p/2}\xi_A(y)^{1-p} dy.
\end{eqnarray*}
This yields that
\begin{equation}
\int_{T_A\cup T_B}\left| \frac{\partial \Pi \vec{u}}{\partial \vec{x}} \right|^p
d\vec{x}=4\tau \sin{\frac{\pi}{2N}}\int_0^1
\Upsilon(y)^{p/2}\xi_A(y)^{1-p}dy(1+\iota_1). \label{TABenergy}
\end{equation}
Similar calculations yield
$$
\Big|\frac{\partial \Pi \vec{u}}{\partial \vec{x}}\Big|^p=
\frac{(2\tau \sin{\frac{\pi}{2N}})^p}{(\det \nabla
\vec{x})^p}\Upsilon(\hat{x}_1)^{p/2}(1+O(\tau^2+N^{-2})),
$$
and
$$
\Big|\frac{\partial \Pi \vec{u}}{\partial \vec{x}}\Big|^p=
\frac{(2\tau \sin{\frac{\pi}{2N}})^p}{(\det \nabla
\vec{x})^p}\Upsilon(\hat{x}_2)^{p/2}(1+O(\tau^2+N^{-2})),$$
on the elements of types C and D respectively.
Thus it follows from Lemma~\ref{wgp} or more precisely  \eqref{detC} and
\eqref{detD} that
\begin{eqnarray*}
&&\int_{T_C\cup T_D}\left| \frac{\partial  \Pi\vec{u}}{\partial \vec{x}} \right|^p
d\vec{x}\\ &=&2\tau \sin{\frac{\pi}{2N}}\int_{\hat{T}}
\big(\Upsilon(\hat{x}_1)^{p/2}\xi_A(\hat{x}_1)^{1-p}+\Upsilon(\hat{x}_2)^{p/2}
\xi_A(\hat{x}_2)^{1-p}\big)d\hat{\vec{x}}(1+O(\tau^2+N^{-2}))\\
&=&4\tau \sin{\frac{\pi}{2N}}\int_0^1 (1-y)\Upsilon(y)^{p/2}\xi_A(y)^{1-p}
dy(1+O(\tau^2+N^{-2})).
\end{eqnarray*}
Comparing this with the corresponding estimate on element $T_B$ (compare \eqref{TBenergy}),
we are led to (compare \eqref{TABenergy})
\begin{equation}
\int_{T_A\cup T_C \cup T_D}\left| \frac{\partial \Pi \vec{u}}{\partial \vec{x}} \right|^p
d\vec{x}=4\tau \sin{\frac{\pi}{2N}}\int_0^1
\Upsilon(y)^{p/2}\xi_A(y)^{1-p}dy(1+\iota_1). \label{TACDenergy}
\end{equation}

Similarly, by Lemma~\ref{wgp}, one has
\begin{eqnarray*}
\int_{T_A}\left| \frac{\partial \vec{u}}{\partial \vec{x}} \right|^p d\vec{x}
&=&\int_{T_A}\Big(r^2(|\vec{x}|)+r'(|\vec{x}|)^2
|\vec{x}|^2\Big)^{p/2}|\vec{x}|^{-p}d \vec{x}\\
&=&\int_{T}\Upsilon(\hat{x}_1+\hat{x}_2)^{p/2}\xi_A(\hat{x}_1+\hat{x}_2)^{-p}\det \nabla \vec{x}d\hat{\vec{x}}(1+O(N^{-2}))\\
&=&4\tau \sin{\frac{\pi}{2N}}\int_0^1 y\Upsilon(y)^{p/2}\xi_A(y)^{-p}(\xi_A(y)+
2\tau \sin^2\frac{\pi}{4N})dy(1+\iota_1), \\
\int_{T_B}\left| \frac{\partial  \vec{u}}{\partial \vec{x}} \right|^p d\vec{x}
&=&4\tau \sin{\frac{\pi}{2N}}\int_0^1 y\Upsilon(1-y)^{p/2}\xi_A(1-y)^{1-p}dy(1+\iota_1)\\
&=&4\tau \sin{\frac{\pi}{2N}}\int_0^1 (1-y)\Upsilon(y)^{p/2}\xi_A(y)^{1-p} dy(1+\iota_1),
\end{eqnarray*}
which yield
\begin{eqnarray} \label{E1AB}
\int_{T_A\cup T_B}\left| \frac{\partial \vec{u}}{\partial \vec{x}} \right|^p d\vec{x}
&=&4\tau \sin{\frac{\pi}{2N}}\int_0^1\Upsilon(y)^{p/2}\xi_A(y)^{1-p} dy(1+\iota_1)
\nonumber \\
&&+4\tau \sin{\frac{\pi}{2N}}\int_0^1y\Upsilon(y)^{p/2}\xi_A(y)^{-p}
2\tau \sin^2\frac{\pi}{4N}dy \nonumber \\
&=&4\tau \sin{\frac{\pi}{2N}}\int_0^1\Upsilon(y)^{p/2}\xi_A(y)^{1-p} dy(1+\iota_1),
\end{eqnarray}
or, since $T_B = T_C \cup T_D$, equivalently
\begin{equation}\label{E1ACD}
\int_{T_A\cup T_C \cup T_D}\left| \frac{\partial \vec{u}}{\partial \vec{x}} \right|^p
d\vec{x}=4\tau \sin{\frac{\pi}{2N}}\int_0^1\Upsilon(y)^{p/2}\xi_A(y)^{1-p} dy(1+\iota_1).
\end{equation}
Thus \eqref{error E13} follows as a direct consequence of \eqref{TABenergy}-\eqref{E1ACD}.
Moreover,
\begin{eqnarray*}
E_1(\vec{u};\Omega^{\mathcal{J}}_{(\epsilon,\tau)})&=&
4\omega \tau  N\sin{\frac{\pi}{2N}}\int_0^1\Upsilon(y)^{p/2}\xi_A(y)^{1-p} dy(1+\iota_1)\\
&=&2\pi \omega \tau \int_0^1\Upsilon(y)^{p/2}\xi_A(y)^{1-p} dy(1+\iota_1)\\
&=&2 \pi \omega \int_\epsilon^{\epsilon+\tau}
(r^2(t)+r'(t)^2t^2)^{p/2} t^{1-p} dt(1+\iota_1)\\
&=&E_1(\vec{u};\Omega_{(\epsilon,\tau)})(1+\iota_1)\\
&\sim &A(\epsilon,\tau).
\end{eqnarray*}
Thus we get
$$
E_1(\vec{u};\Omega^{\mathcal{J}}_{(\epsilon,\tau)})=
E_1(\vec{u};\Omega_{(\epsilon,\tau)})(1+O(\tau^2+N^{-2}+\frac{\epsilon}{\tau}N^{-4}))
\sim A(\epsilon,\tau).
$$

On the other hand, by \eqref{det hPiu} and
$g(\cdot) \in C^2(0,+\infty)$, on the elements of types A, B, C and D, we have
$$
g\Big(\det \frac{\partial \Pi \vec{u}}{\partial \vec{x}}\Big)=
g(\det \frac{\partial \vec{u}}{\partial \vec{x}})+\frac{N^{-1}\iota_2}{
\det \nabla \vec{x}},$$
where $\iota_2=O(\tau^3+(\epsilon+\tau)\tau N^{-2} +N^{-4})$.
Hence
$$
|E_2(\Pi \vec{u};\Omega^{\mathcal{J}}_{(\epsilon,\tau)})-
E_2(\vec{u};\Omega^{\mathcal{J}}_{(\epsilon,\tau)})|=
N\int_{T_A\cup T_B}\frac{1}{\det \nabla \vec{x}}d\vec{x}N^{-1}\iota_2
=\iota_2,
$$
which is \eqref{error E23}. Similar arguments and the fact that $g(\det \nabla \vec{u})
\ge g_0>0$ lead to
$$
E_2(\vec{u};\Omega^{\mathcal{J}}_{(\epsilon,\tau)})=
E_2(\vec{u};\Omega_{(\epsilon,\tau)})(1+\iota_1)\sim \epsilon\tau+\tau^2,
$$
hence \eqref{Eo} holds, i.e., $E(\vec{u};\Omega^{\mathcal{J}}_{(\epsilon,\tau)})
\sim E(\vec{u};\Omega_{(\epsilon,\tau)})\sim A(\epsilon,\tau)
\sim \max\{\epsilon,\tau\}^{1-p}\tau$.
Thus, one has
$$
\frac{|E_2(\Pi \vec{u};\Omega^{\mathcal{J}}_{(\epsilon,\tau)})-E_2(\vec{u};
\Omega^{\mathcal{J}}_{(\epsilon,\tau)})|}{
E(\vec{u};\Omega^{\mathcal{J}}_{(\epsilon,\tau)})}=\frac{\iota_2}{A(\epsilon,\tau)}
\sim \max\{\epsilon,\tau\}^{p-1}\iota_2/\tau,
$$
which together with \eqref{error E13} gives \eqref{error E3}.
\hfill $\square$

\section{Meshing strategy and convergence theorem}

In this section, by applying the error estimates given in \S~3 and the
orientation-preservation condition given in \cite{detp} (see also
Remark~\ref{rm-orientation}), we first establish a meshing
strategy on the domain $\Omega_{\epsilon_0}=B_1(\vec{0})\setminus B_{\epsilon_0}(\vec{0})$.
The mesh is assumed to be introduced by $\mathcal{J}=\bigcup \limits_{i=0}^{m}
\Omega^{\mathcal{J}}_{(\epsilon_i,\tau_i)}$ with each $\Omega^{\mathcal{J}}_{(\epsilon_i,\tau_i)}$ being either a radially symmetric
mesh on a circular ring domain $\{x : \epsilon_i \leq |x| \leq \epsilon_{i+1}\}$
consisting of $2N_i$ iso-parametric quadratic finite elements of types A and B
as shown in Figure~\ref{tuli}, or a slightly modified mesh consisting of $3N_i$
iso-parametric quadratic finite elements of types A, C and D as shown in
Figure~\ref{tuli1}. Our purpose is,
given a far field reference mesh size $h>0$, to find $m$, and $\epsilon_i$,
$N_i$, $i=0,1,\cdots, m$, so that (1): the finite element interpolation function of the
cavity deformation is orientation preserving; (2): the relative errors of the
elastic energy $E$ on the circular ring domains
$\{x : \epsilon_i \leq |x| \leq \epsilon_{i+1}\}$ are all of the
order $O(h^{2})$; and (3): at the same time, the absolute errors of the elastic
energy on each of the circular domains are of the same order. The last
requirement, which can be realized by making $(\epsilon_i + \tau_i)^{2-p}-
\epsilon_i^{2-p} \sim E_1(\vec{u};\Omega_{(\epsilon_i,\tau_i)}) \sim
E(\vec{u};\Omega_{(\epsilon_i,\tau_i)})$ (see \eqref{Eo})
to be the same order, implies that the absolute error of the
elastic energy is in some sense equi-distributed in the radial direction.

Let $\epsilon_0< \epsilon_1 \cdots < \epsilon_i < \cdots< \epsilon_{m}<
\epsilon_{m+1}=1.0$, let $\tau_{i}= \epsilon_{i+1} -\epsilon_i$,
and let $N_i$ be the number of the nodes on both the inner and outer boundaries
of the circular ring domain $\Omega_{(\epsilon_i,\tau_i)}$ introduced by
$\Omega^{\mathcal{J}}_{(\epsilon_i,\tau_i)}$ (see
Figure~\ref{tuli}). For the simplicity of the finite element coding, we
require that either $N_i=2N_{i+1}$ or $N_i=N_{i+1}$. By Theorem 3.7 of \cite{detp}
(see also Remark~\ref{rm-orientation}), to preserve the orientation of the
finite element interpolation
functions, $\epsilon_i$, $\tau_i$, $N_i$ must satisfy the conditions
$\tau_i\le C_1\epsilon_i^{1/2}$, and $N_i^{-1} \leq C_2(\epsilon_i
\tau_i)^{1/4}$, where $C_1$, $C_2$ are constants depending on the solution
$r(R)$. On the other hand, by \eqref{error E3}, to ensure the relative error of
$E(\Pi \vec{u};\Omega^\mathcal{J}_{(\epsilon_i,\tau_i)})$ to be the order of $O(h^2)$,
a necessary condition is that $\tau_i = O(h)$, $\forall i$. Moreover, since
$A(\epsilon_m, \tau_m)=1-(1-\tau_m)^{2-p}\le (2-p)2^{p-1}\tau_m=O(h)$ if
$\epsilon_m>\frac{1}{2}$,
it is natural to require $A(\epsilon_i, \tau_i)\le Ch$, for all $0 \le i \le m$,
which imposes a condition on the layer's thickness: $\tau_i \le d(\epsilon_i,h)$,
where $d(x,h):=(x^{2-p}+Ch)^{\frac{1}{2-p}}-x$ is defined by $A(x,d(x,h))=Ch$.
As is shown in \cite{SuLiRectan}, $\tau_i\le d(\epsilon_i,h)$ leads to
$\tau_i=O(h)$ as $h \rightarrow 0$.

For given positive constants $C_1$, $C_2$, $C\ge (2-p)2^{p-1}$, $h\le
\min\{\frac{2-p}{2^{2-p}C}, \frac{2-p}{2^{p-1}C}\}$ (see \cite{SuLiRectan}),
$A_1<A_2$ satisfying $[(A_2h)^{-1},(A_1h)^{-1}]\cap \mathbb{Z}_+ \ne \emptyset$,
the analysis above leads to the following meshing strategy.

\vskip 3mm
{\bf A meshing strategy of $\{\Omega^{\mathcal{J}_h}_{(\epsilon_i,\tau_i)}\}_{i=0}^{m}$:}
\begin{description}
\item[(1)]
Set $\tau_0=\min \{C_1\epsilon_0^{1/2},d(\epsilon_0,h)\}$.
Take $\tilde{N}_m \in [(A_2h)^{-1}, (A_1h)^{-1}] \cap \mathbb{Z}_+ $.
Let $\bar{N}_0= \min \{ N \in \mathbb{Z}_+: N^{-1}\le
\min\{ C_2(\epsilon_0\tau_0)^{1/4},(C^2\tau_0h^2)^{1/4}\} \}$.
Set $k=\min \{ j: 2^j
\tilde{N}_m \ge \bar{N}_0 \}$, and $N_0={2^k}\tilde{N}_m$.
\item[(2)] Set $k_0=0$. For $i \ge 1$, set $\epsilon_i=\epsilon_{i-1}+\tau_{i-1}$, and
\begin{equation}\label{tau}
\tau_i=\min\{1-\epsilon_i,C_1\epsilon_i^{1/2},d(\epsilon_i,h)\}.
\end{equation}
If $\tau_i=1-\epsilon_i$, set $m=i$. The least admissible
$N_i$ such that $N_i^{-1}\le \min\{ C_2(\epsilon_i\tau_i)^{1/4},(C^2\tau_ih^2)^{1/4}\}$
is determined as follows:
\begin{enumerate}
\item If $k_{i-1}<k$, set $\bar{N}_i=\frac{N_{i-1}}{2}$. If
$\bar{N}_i^{-1}\le \min\{ C_2(\epsilon_i\tau_i)^{1/4},(C^2\tau_ih^2)^{1/4}\}$,
then set $k_i=k_{i-1}+1$, $N_i=\bar{N}_i$; otherwise, set $k_i=k_{i-1}$,
$N_i=N_{i-1}$.
\item If $k_{i-1}=k$, set $k_i=k_{i-1}$, $N_i=N_{i-1}$.
\end{enumerate}
\item[(3)] Whenever $N_i=2N_{i+1}$ occurs, each type B element in the circular ring
$\Omega^\mathcal{J}_{(\epsilon_i,\tau_i)}$ is divided into a pair of elements of
types C and D by introducing a straight line right in the middle along the radial
direction (see Figure~\ref{tuli1}).
\end{description}

\begin{remark}
The step (3) above can be viewed as a conforming process.
The quadratic iso-parametric finite element function space established on a mesh produced
according to the above meshing strategy is a conforming finite element function space.
The analysis in \S~3 shows that the orientation-preservation conditions as well as
the error bounds of the interpolation function of the
cavity deformation are not jeopardised by the conforming process.
\end{remark}

On a mesh $\{\Omega^{\mathcal{J}_h}_{(\epsilon_i,\tau_i)}\}_{i=0}^{m}$ produced according
to the above meshing strategy, we have the following results.

\begin{theorem}\label{theorem err=O(h^2)}
Let $\Omega_h=\cup_{i=0}^m \Omega^{\mathcal{J}_h}_{(\epsilon_i, \tau_i)}$,
let $\vec{u}(\vec{x})$ be the radially symmetric cavity deformation on
$\Omega_{\epsilon_0}$, then
$\det \nabla \Pi \vec{u}(\vec{x})>0$ a.e.~on $\Omega_h$. Moreover, the error of the
elastic energy satisfies
\begin{equation}\label{errE<=O(h^2)}
E(\Pi \vec{u};\Omega_h)=\sum\limits_{i=0}^m E(\Pi \vec{u};
\Omega^{\mathcal{J}_h}_{(\epsilon_i, \tau_i)})=E(\vec{u}; \Omega_{\epsilon_0} )(1+O(h^2)).
\end{equation}
\end{theorem}
\textbf{Proof.} Since $\tilde{N}_m\sim 1/h$, $N_i\ge \tilde{N}_m$, it follows that
$N_i^{-1}=O(h)$. It is easily verified that the orientation-preservation conditions
$\tau_i\le C_1\epsilon_i^{1/2}$ and $N_i^{-1}\le C_2(\epsilon_i\tau_i)^{1/4}$ are
satisfied and $\tau_i=O(h), \frac{1}{\tau_i}N_i^{-4}=O(h^2)$. Thus,
$\det \nabla \Pi \vec{u}(\vec{x})>0$
a.e. on $\Omega_h$ by Theorem 3.7 of \cite{detp} (see also Remark~\ref{rm-orientation}).
On the other hand, it follows from \eqref{Eo} and \eqref{error E3} that, for all $i$,
$E(\Pi \vec{u}; \Omega^{\mathcal{J}_h}_{(\epsilon_i,
\tau_i)})=E(\vec{u}; \Omega_{(\epsilon_i, \tau_i)})(1+O(h^2))$, which yield
\eqref{errE<=O(h^2)}.
\hfill $\square$

\begin{theorem}
Let $\Omega_{h_k}=\cup_{i=0}^m \Omega^{\mathcal{J}_{h_k}}_{(\epsilon_i, \tau_i)}$ with
$\lim\limits_{k \rightarrow \infty}h_k = 0$, let $\mathcal{A}_{h_k}$ be the corresponding
conforming finite element function spaces consist of piecewise quadratic iso-parametric
functions satisfying the boundary condition $\vec{v}_{h_k}(\vec{x})=\lambda \vec{x}$ for
all the mesh nodes on the boundary $\Gamma_0$, and let $\vec{u}_{h_k}$ be a minimizer
of $E(\cdot)$ in $\mathcal{A}_{h_k}$. Suppose that the radially symmetric cavity
deformation $\vec{u}$
is the unique minimizer of $E(\cdot)$ in $\mathcal{A}$, then
$\vec{u}_{h_k}\chi_{\Omega_{\epsilon_0}} \rightarrow \vec{u}$ in
$W^{1,p}(\Omega_{\epsilon_0},\mathbb{R}^2)$.
\end{theorem}
\textbf{Proof.} By Theorem~\ref{theorem err=O(h^2)} and the assumption that
$\vec{u}$ is the energy minimizer, we conclude that
\begin{equation}\label{FEMenergyconverges}
\lim_{k \rightarrow \infty} E(\vec{u}_{h_k}) = E(\vec{u}) =
\inf_{\vec{v}\in \mathcal{A}} E(\vec{v}),
\end{equation}
and in particular $\{\|\nabla \vec{u}_{h_k}\|_p\}$ is bounded, since $g> 0$.
This, by the boundary condition and the Poincar\'{e} inequality \cite{Morrey},
implies that $\{\vec{u}_{h_k}\chi_{\Omega_{\epsilon_0}}\}$ is bounded in
$W^{1,p}(\Omega_{\epsilon_0},\mathbb{R}^{2})$. Consequently, by the
De La Vall\'{e}e Poussin
criterion (\cite{Meyer}), both $\{\nabla (\vec{u}_{h_k}\chi_{\Omega_{\epsilon_0}})\}$ and
$\{\det \nabla (\vec{u}_{h_k}\chi_{\Omega_{\epsilon_0}})\}$ are equi-integrable,
since $p>1$ and $g$ is a convex function satisfying \eqref{asmp g}.
Thus, there exist a subsequence $\vec{u}_{h_k}\chi_{\Omega_{\epsilon_0}}$ (not relabelled),
$\vec{\tilde{u}} \in W^{1,p}(\Omega_{\epsilon_0},\mathbb{R}^{2})$ and
$\vartheta \in L^1(\Omega_{\epsilon_0})$ such that
$$
\vec{u}_{h_k}\chi_{\Omega_{\epsilon_0}} \rightharpoonup \vec{\tilde{u}} \in
W^{1,p}(\Omega_{\epsilon_0},\mathbb{R}^2), \; \vec{u}_{h_k}\chi_{\Omega_{\epsilon_0}}
\rightarrow \vec{\tilde{u}}~a.e., \; \det \nabla (\vec{u}_{h_k}\chi_{\Omega_{\epsilon_0}})
\rightharpoonup \vartheta \in L^1(\Omega_{\epsilon_0}).
$$

Clearly $\vartheta \ge 0$ a.e., we claim that $\vartheta >0$ a.e..
Suppose otherwise, i.e. if $\vartheta$ were zero in a
set $A$ of positive measure, then one would have $\int_A |\det \nabla
\vec{u}_{h_k}| \rightarrow 0$ and $\det \nabla \vec{u}_{h_k} \rightarrow 0$ a.e.
in $A$. Hence, by the assumption of $g$, one would have
$g(\det \nabla \vec{u}_{h_k})\rightarrow \infty$ a.e. in $A$, and as a consequence
$E(\vec{u}_{h_k}) \rightarrow \infty$, which is a contradiction.

The fact that $\vec{u}_{h_k} \chi_{\Omega_{\epsilon_0}}\rightharpoonup \vec{\tilde{u}}$
in $W^{1,p}(\Omega_{\epsilon_0},\mathbb{R}^2)$ implies that $\nabla \vec{u}_{h_k}
\chi_{\Omega_{\epsilon_0}} \rightharpoonup \nabla \vec{\tilde{u}}$ in
$L^p(\Omega_{\epsilon_0}, \mathbb{R}^{2\times 2})$, $\cof \nabla \vec{u}_{h_k}
\chi_{\Omega_{\epsilon_0}} \rightharpoonup \cof \nabla \vec{\tilde{u}}$ in
$L^1(\Omega_{\epsilon_0}, \mathbb{R}^{2\times 2})$. In addition,
since $\vec{u}_{h_k}$ is continuous, by Theorem 3 of \cite{Henao2011},
$\mathcal{E}(\vec{u}_{h_k})=0$, where
$$
\mathcal{E}(\vec{v}):=\sup \{ \mathcal{E}(\vec{v},\vec{f}):\vec{f} \in
C_c^{\infty}(\Omega \times \mathbb{R}^n,\mathbb{R}^n), ~||\vec{f}||_{\infty} \le 1\},
$$
$$
\mathcal{E}(\vec{v},\vec{f}):=\int_\Omega [\cof \nabla \vec{v}
(\vec{x}) \cdot \nabla_{\vec{x}} \vec{f}(\vec{x},\vec{v}(\vec{x}))+
\det \nabla \vec{v}(\vec{x}) \sdiv_{\vec{v}} \vec{f}(\vec{x},\vec{v}(\vec{x}))]d\vec{x},
$$
and where $\nabla_{\vec{x}}$ and $\sdiv_{\vec{v}}$ denote the gradient and divergence
of $f(\vec{x},\vec{v})$ with respect to $\vec{x}$ and $\vec{v}$ respectively.
Thus, by Theorem 2 of \cite{Henao2010} $\vec{\tilde{u}}$ is one-to-one almost everywhere,
and by Theorem 3 of \cite{Henao2010} $\vartheta =\det \nabla \vec{\tilde{u}}$, a.e.~and
$\mathcal{E}(\vec{\tilde{u}})=0$. Moreover,
by the boundary condition of $\vec{u}_{h_k}$ on $\Gamma_0$ and the relationship between
$\Omega_{h_k}$ and $\Omega_{\epsilon_0}$ we conclude that
$\vec{\tilde{u}}|_{\Gamma_0}=\lim_{k \rightarrow \infty}\vec{u}_{h_k}|_{\Gamma_0}=
\lambda \vec{x}$,
hence $\vec{\tilde{u}} \in \mathcal {A}$.
Thus, by the lower semi-continuity theorem (\cite{Ball81}, Theorem 5.4) and
\eqref{FEMenergyconverges}, we obtain that $\inf\limits_{\vec{v}\in \mathcal{A}} E(\vec{v}) \le
E(\vec{\tilde{u}})\le \liminf\limits_{k \rightarrow \infty}E(\vec{u}_{h_k})=
\inf\limits_{\vec{v}\in \mathcal{A}} E(\vec{v})$, which implies that $\vec{\tilde{u}}=\vec{u}$
is the unique minimizer of $E(\cdot)$ in $\mathcal{A}$.

On the other hand, it follows from the convexity of $g$ that
\begin{eqnarray*}
E(\vec{u})-\omega\int_{\Omega_{\epsilon_0}}|\nabla \vec{u}|^p dx
dx&=&\int_{\Omega_{\epsilon_0}}g(\det \nabla \vec{u})dx\\ &\le&\liminf_{k
\rightarrow \infty}\int_{\Omega_{\epsilon_0}}g(\det \nabla \vec{u}_{h_k})dx\\
&=&\liminf_{k \rightarrow \infty}(E(\vec{u}_{h_k})-\omega\int_{\Omega_{\epsilon_0}}|\nabla
\vec{u}_{h_k}|^pdx)\\
&=&E(\vec{u})-\omega\limsup_{k \rightarrow \infty}\int_{\Omega_{\epsilon_0}}|\nabla
\vec{u}_{h_k}|^pdx.
\end{eqnarray*}
This implies that $\|\nabla \vec{u}\|_p=\lim\limits_{k \rightarrow \infty}
\|\nabla \vec{u}_{h_k}\|_p$, which
together with $\vec{u}_{h_k} \chi_{\Omega_{\epsilon_0}}\rightharpoonup \vec{u}\in
W^{1,p}(\Omega_{\epsilon_0},\mathbb{R}^2)$ yields $\vec{u}_{h_k}\chi_{\Omega_{\epsilon_0}}
\rightarrow \vec{u} \in W^{1,p}(\Omega_{\epsilon_0},\mathbb{R}^2)$ (see \cite{Evans}).
\hfill $\square$

\section{Numerical experiments and results}

In this section, the numerical results are presented to illustrate the efficiency of our
meshing strategy. Before proceeding, we notice that, in the meshing strategy, there are
two solution-dependent constants $C_1$ and $C_2$, which are not known a priori. However,
in applications, we can always start with $C_1:= d(\epsilon_0,h)
\epsilon_0^{-1/2}$ and $C_2:=\tilde{N}_m^{-1} (\epsilon_0d(\epsilon_0,h))^{-1/4}$,
which are the least $C_1$ and $C_2$ such that the orientation-preservation conditions
will practically not affect the mesh produced. The numerical solutions on an improper
mesh with the constants $C_1$ or $C_2$ too large might still capture the cavitation
phenomenon, but would typically fail to be orientation preserving, which
can be most easily detected on the corners of the elements on the inner boundary of
$\Omega_{\epsilon_0}$. Similar as in \cite{SuLiRectan}, whenever the failure of
orientation preservation is detected, the constant $C_1$ or $C_2$ or both should be
reduced, say by half, or simply increase $N_0$ instead, say by doubling, and the
process, repeat if necessary, will efficiently produce a proper mesh in the end.

The energy density in our numerical experiments is given by
\eqref{energy density} with $p=3/2$, $\omega=2/3$, and
$g(x)=2^{-1/4}(\frac12(x-1)^2+\frac{1}{x})$, the domain is
$\Omega_{\epsilon_0} = B_1(\vec{0}) \setminus B_{\epsilon_0}(\vec{0}) \subseteq \mathbb{R}^2$
with a displacement boundary condition $\vec{u}_0(\vec{x})=2\vec{x}$ given on
$\Gamma_0=\partial B_1(\vec{0})$ and a traction free boundary condition given on
$\Gamma_1=\{\vec{x}: |\vec{x}|=\epsilon_0\}$,
and the meshes used are shown in Table~\ref{e-2} and Table~\ref{e-4}, which are produced by the
meshing strategy with $C=2$, $C_1=0.9$, $C_2=0.5$, $A_1=0.8$, $A_2=1$ for $\epsilon_0=0.01$,
$\epsilon_0=0.0001$ and various $h$.

\begin{table}[ht!]
\small
\begin{minipage}[l]{0.475\textwidth}
\caption{$\epsilon_0=0.01$.}\label{e-2}
\vspace*{1mm}\hspace*{1mm}
  \begin{tabular}{|l|l|l|l|l|l|l|l|l|l|l|}
\hline
$h$&$\min \tau_i$&$\max \tau_i$&$m$& $N_h$\\\hline
0.06&0.0384&0.2112&7&15\\\hline
0.04&0.0224&0.1504&11&20\\\hline
0.03 &0.0156&0.1164&14&27\\\hline
0.02&0.0096&0.0768&22&40\\\hline
0.01&0.0044&0.0396&44&80\\\hline
\end{tabular}
\end{minipage}\hspace*{5mm}
	\begin{minipage}[l]{0.525\textwidth}
\caption{$\epsilon_0=0.0001$.}\label{e-4}
\vspace*{1mm}
\begin{tabular}{|l|l|l|l|l|l|l|l|l|l|l|}
\hline
$h$&$\min \tau_i$&$\max \tau_i$&$m$&$N_0$&$N_m$\\\hline
0.06&0.009&0.21&8&16&64 \\\hline
0.04&0.008&0.1488&12&20&80\\\hline
0.03&0.0048&0.1128&16 &27&108 \\\hline
0.02&0.0024&0.076&24&46& 92 \\\hline
0.01&0.0008&0.0392&49&80 &160\\\hline
\end{tabular}
\end{minipage}
\end{table}

\begin{figure}[ht!]
	\begin{minipage}[l]{0.5\textwidth}
    \centering    \includegraphics[width=2.65in]{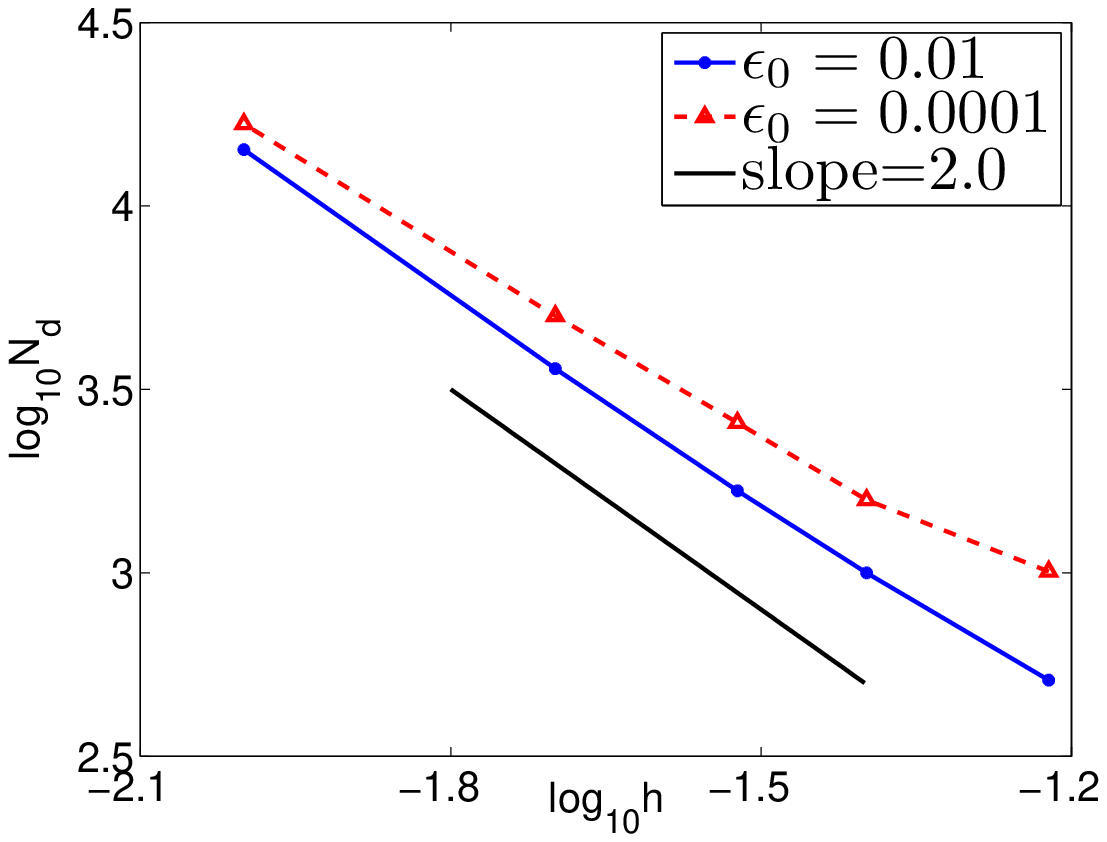}
	\vspace*{-4mm}
	\caption{$N_d \sim h^{-2}$.}\label{Nd_h}
\end{minipage}\hspace*{2mm}
\begin{minipage}[l]{0.5\textwidth}
\centering   \includegraphics[width=2.65in]{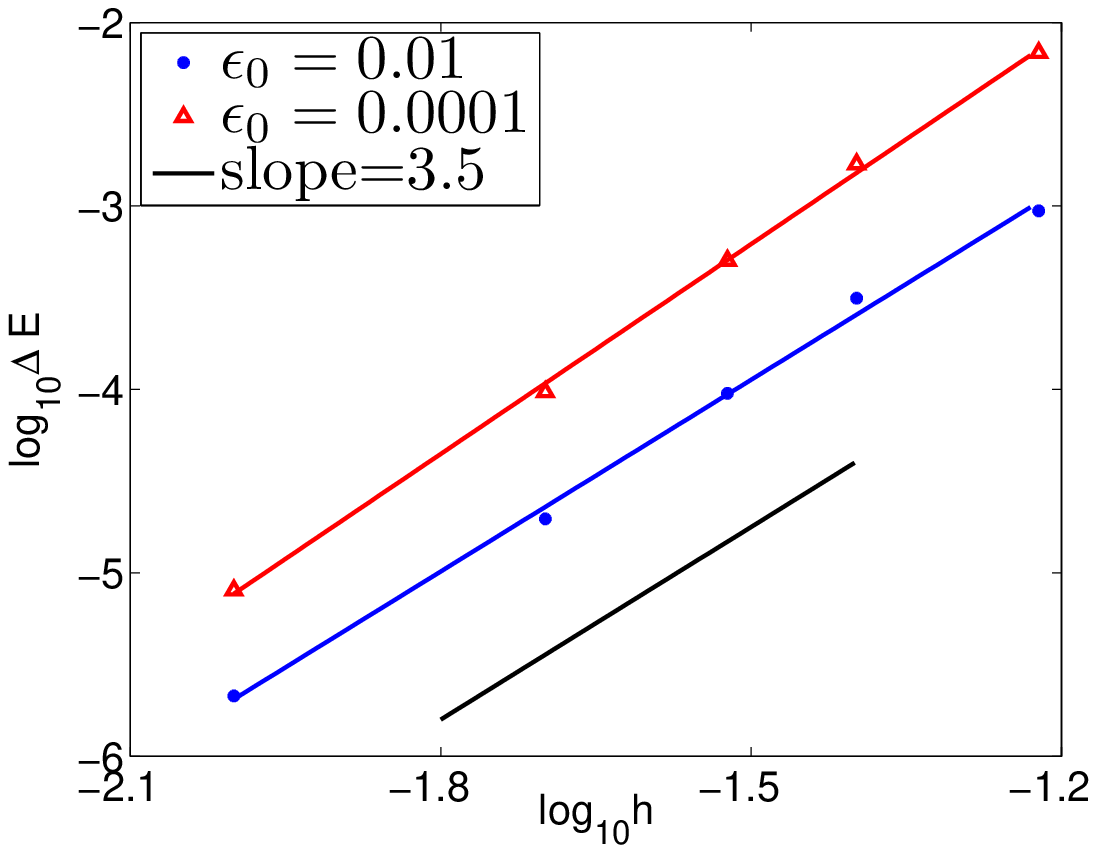}
\vspace*{-4mm}
\caption{The energy error.}\label{energy error}
\end{minipage}
\end{figure}

It happens that, for $\epsilon_0=0.01$, $N_i=N_h$ on each of the $m+1$ mesh layers,
while for $\epsilon_0=0.0001$, the $N_i=2N_{i+1}$ do occur in several of the innermost layers,
however in both cases the total degrees of freedom $N_d$ is asymptotically a quadratic
function of $h^{-1}$ as shown in Figure \ref{Nd_h}.

\begin{figure}[ht!]
	\begin{minipage}[l]{0.5\textwidth}
    \centering    \includegraphics[width=2.65in]{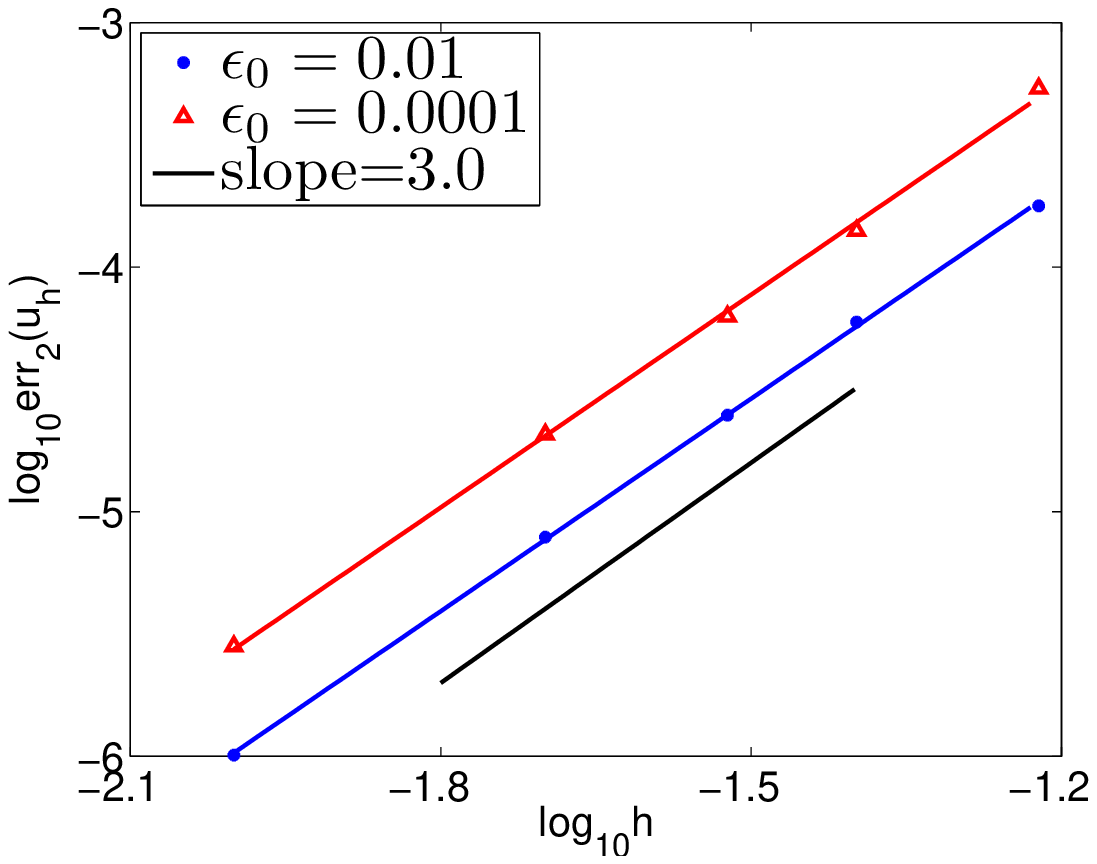}
	\vspace*{-4mm}
	\caption{The $L^2$ error of $\vec{u}_h$.}\label{l2 error}
\end{minipage}\hspace*{2mm}
\begin{minipage}[l]{0.5\textwidth}
\centering   \includegraphics[width=2.65in]{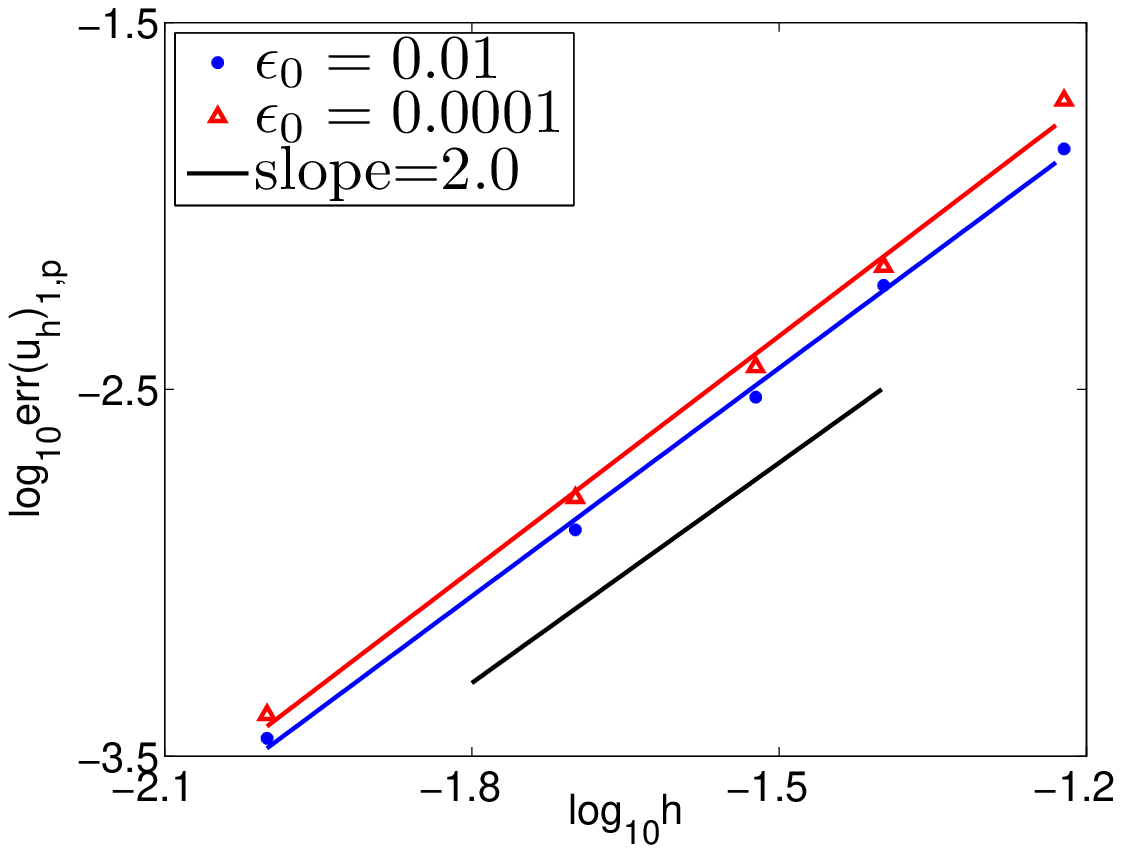}
\vspace*{-4mm}
\caption{The $W^{1,p}$ error of $\vec{u}_h$.}\label{w1p error}
\end{minipage}
\end{figure}

The convergence behavior of the elastic energy is shown in Figure
\ref{energy error}, where it is clearly seen that the convergence rate of the elastic
energy of the finite element solutions $\vec{u}_h$ is more than one order higher than that
of one could standardly expect from a quadratic approximation (see also \eqref{error E3})
showing that the method probably has some kind of super-convergence potential.
In Figure~\ref{l2 error} and Figure~\ref{w1p error}, we see that
$\|\vec{u}-\vec{u}_h\|_{0,2}=O(h^3)$ and $\|\vec{u}-\vec{u}_h\|_{1,p}=O(h^2)$
respectively, which show that our meshing strategy is optimal in the sense that
the optimal order of convergence rates in $\|\cdot\|_{0,2}$ and $\|\cdot \|_{1,p}$ norms
can be achieved with the quadratic iso-parametric FEM, recalling that $N_d \sim h^{-2}$.

\begin{figure}[h!]
     \begin{minipage}[c]{0.55\textwidth}\hspace{-2mm}
     \subfigure[Experience based vs. optimized meshes.]{
     \includegraphics[width=6.75cm, height=5cm]{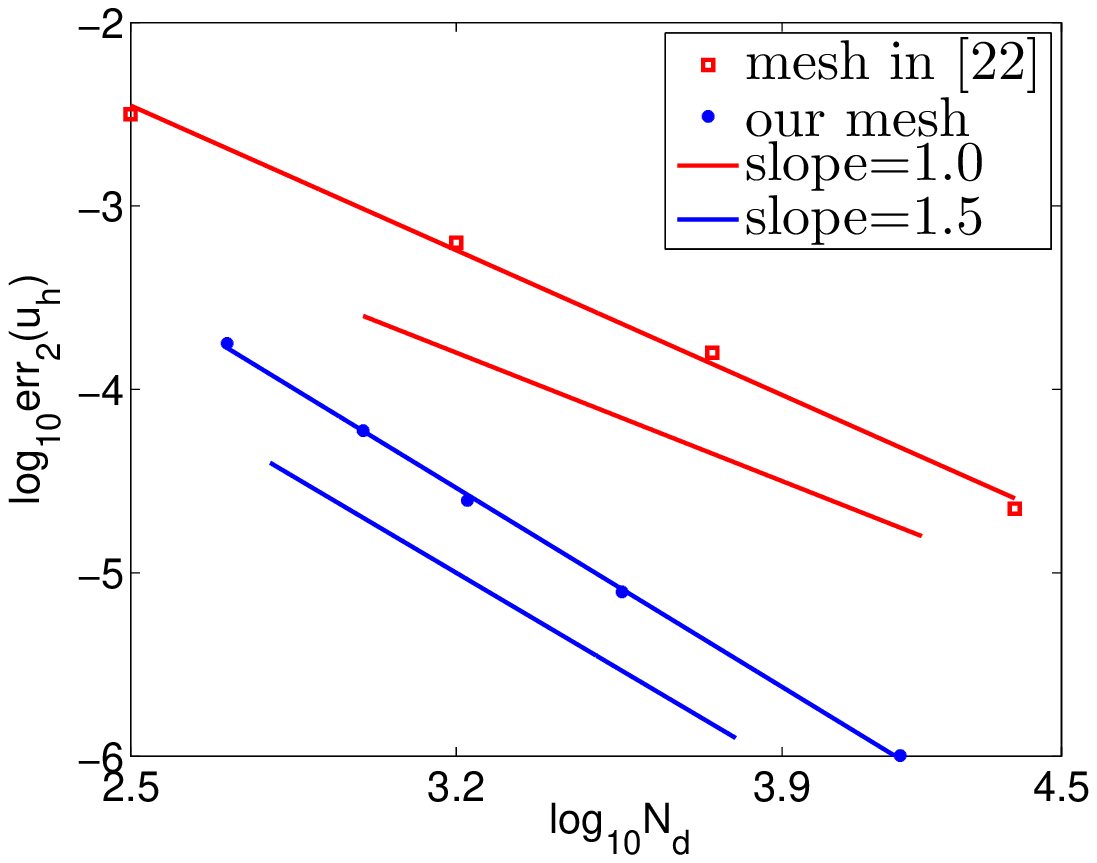}\label{bijiao}}
     \end{minipage}
     \begin{minipage}[c]{0.45\textwidth}\hspace{-8mm}
     \subfigure[Iso- and dual-parametric FEM.]{
     \includegraphics[width=6.75cm, height=5cm]{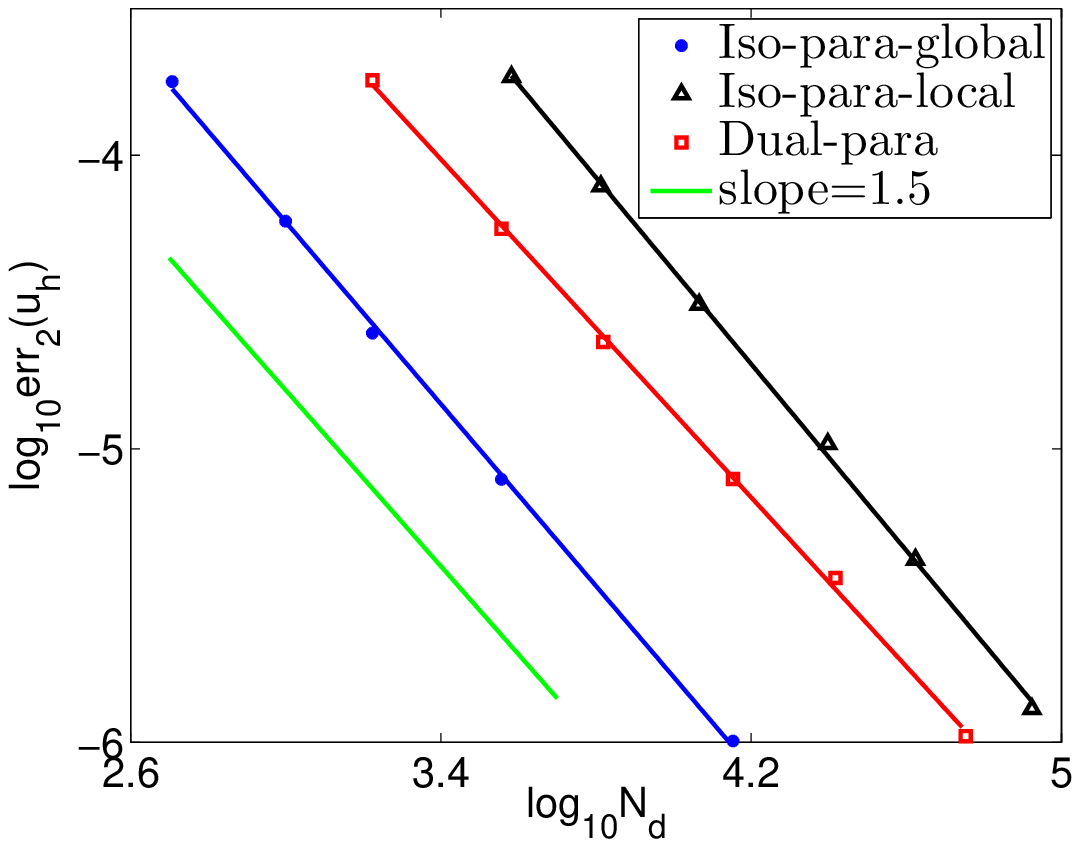}\label{l2_34}}
     \end{minipage}
     \vspace*{-4mm}
 \caption{The comparison of the $L^2$ errors of numerical solutions ($\epsilon_0=0.01$).}
 \label{bijiao-2}
\end{figure}

Figure~\ref{bijiao} compares the $L^2$ error of the numerical
cavity solutions obtained on the meshes produced by our meshing strategy and on
the meshes provided according to the limited numerical experiences given in
\cite{Lian and Li iso}. Figure~\ref{l2_34} displays the $L^2$ error of the numerical
cavity solutions obtained by the quadratic iso-parametric FEM on the globally optimized
meshes and locally optimized ones, where the mesh is optimized by using our meshing
strategy only on $\{\vec{x}: \epsilon_0\le |\vec{x}|\le 0.1\}$, and by the bi-quadratic dual-parametric
FEM on the globally optimized meshes (see \cite{SuLiRectan}), denoted in the legend
as Iso-para-global, Iso-para-local and Dual-para respectively. It is clearly seen
that our meshing strategy efficiently works.

\section{Concluding remarks}

The error estimates obtained in this paper on the quadratic iso-parametric finite
element interpolation functions of the cavity deformations enable us to establish
the meshing strategy in a neighborhood of a pre-existing defect, and consequently
to bound the error of the elastic energy of the conforming finite element cavity
solutions in the order of $O(h^2)$, where $h$ is the
far field mesh size, and further to prove the convergence of the finite element
solutions.

Our numerical experiments show that the convergence behavior of the finite element
solutions in $W^{1,p}$ and $L^2$ norms with respect to $h$ and $N_d$ is
essentially asymptotically independent of $\epsilon_0$.
In fact, for $\epsilon_0=10^{-2}$ and $10^{-4}$, the errors in $W^{1,p}$ and
$L^2$ norms drop to the levels below $10^{-3}$ and $10^{-5}$
respectively when $h=0.01$ and $N_d$ reaches about $1.6 \times 10^4$.
Furthermore, the numerical experiments show that the rate of the elastic energy
error of the numerical cavity solutions reaches the level of $O(h^{3.5})$, indicating that
the numerical solutions obtained by the quadratic iso-parametric FEM on
the meshes produced according to our meshing strategy might
have certain super-convergence character, which yet remains to be explored.
The results suggest that the quadratic iso-parametric
finite element method coupled with our meshing strategy could be considered as a
reliable and efficient tool to compute the cavitation problems in nonlinear elasticity.

\end{document}